\documentclass[14pt]{article}
\usepackage{mathrsfs}
\usepackage{amsthm}
\usepackage{amsmath}
\usepackage{graphicx}
\usepackage{color}
\usepackage{amsfonts}
\newtheorem{theorem}{Theorem}[section]

\newtheorem{lemma}[theorem]{Lemma}
\newtheorem{proposition}[theorem]{Proposition}

\newtheorem{corollary}[theorem]{Corollary}
\numberwithin{equation}{section}
\normalsize

\begin{document}
\title{\textbf{Two limit theorems for the high-dimensional two-stage contact process}}

\author{Xiaofeng Xue \thanks{\textbf{E-mail}: xfxue@bjtu.edu.cn \textbf{Address}: School of Science, Beijing Jiaotong University, Beijing 100044, China.}\\ Beijing Jiaotong University}

\date{}
\maketitle

\noindent {\bf Abstract:}
In this paper we are concerned with the two-stage contact process introduced in \cite{Krone1999} on a high-dimensional lattice. By comparing this process with an auxiliary model which is a linear system, we obtain two limit theorems for this process as the dimension of the lattice grows to infinity. The first theorem is about the upper invariant measure of the process. The second theorem is about asymptotic behavior of the critical value of the process. These two theorems can be considered as extensions of their counterparts for the basic contact processes proved in \cite{Grif1983} and \cite{Schonmann1986}.


\section{Introduction}\label{section one}
In this paper we are concerned with the two-stage contact process on $\mathbb{Z}^d$ introduced in \cite{Krone1999}. First we introduce some notations and definitions for later use. For each
\[
x=\big(x(1),\ldots,x(d)\big)\in \mathbb{Z}^d,
\]
we use $\|x\|$ to denote the $l_1$-norm of $x$, i.e., $\|x\|=\sum_{i=1}^d|x(i)|$. For any $x,y\in \mathbb{Z}^d$, we write $x\sim y$ when and only when $\|x-y\|=1$. In other words, $x\sim y$ means that $x$ and $y$ are neighbors on $\mathbb{Z}^d$. For $1\leq i\leq d$, we use $e_i$ to denote the $i$th elementary unit vector of $\mathbb{Z}^d$, i.e.,
\[
e_i=(0,\ldots,0,\mathop 1\limits_{i \text{th}},0,\ldots,0).
\]
We use $O$ to denote the origin of $\mathbb{Z}^d$, i.e., $O=(0,0,\ldots,0)$.

The two-stage contact process $\{\eta_t\}_{t\geq 0}$ on $\mathbb{Z}^d$ is a continuous-time Markov process with state space $X=\{0,1,2\}^{\mathbb{Z}^d}$ and generator function $\Omega$ given by
\begin{equation}\label{equ 1.1 generator}
\Omega f(\eta)=\sum_{x\in
\mathbb{Z}^d}\sum_{i=0,1,2}H(x,i,\eta)\big[f(\eta^{x,i})-f(\eta)\big]
\end{equation}
for any $\eta\in \{0,1,2\}^{\mathbb{Z}^d}$ and $f\in C(X)$, where
\[
\eta^{x,i}(y)=
\begin{cases}
\eta(y) & \text{~if~}y\neq x,\\
i & \text{~if~}y=x
\end{cases}
\]
and
\[
H(x,i,\eta)=
\begin{cases}
1 & \text{~if~}\eta(x)=2 \text{~and~}i=0,\\
1+\delta & \text{~if~}\eta(x)=1 \text{~and~}i=0,\\
\gamma & \text{~if~}\eta(x)=1 \text{~and~}i=2,\\
\lambda\sum_{y:y\sim x}1_{\{\eta(y)=2\}} & \text{~if~}\eta(x)=0 \text{~and~}i=1,\\
0 & \text{~else}
\end{cases}
\]
for any $x\in \mathbb{Z}^d$ and $i\in \{0,1,2\}$, where $\lambda, \delta, \gamma$ are positive constants and $1_A$ is the indicator function of the event $A$ that $1_A=1$ on the event $A$ while $1_A=0$ on the complementary set of $A$.

Intuitively, the two-stage contact process describes the spread of an epidemic on the graph $\mathbb{Z}^d$. The vertices in state $0$ are healthy and vertices in state $1$ are semi-infected while vertices in state $2$ are fully-infected. A fully-infected vertex waits for an exponential time with rate $1$ to become healthy. A semi-infected vertex waits for an exponential time with rate $1+\delta$ to become healthy while waits for an exponential time with rate $\gamma$ to become fully-infected, depending on which moment comes first.  A healthy vertex is infected to become semi-infected at rate proportional to the number of fully-infected neighbors.

The two-stage contact process $\{\eta_t\}_{t\geq 0}$ is introduced in \cite{Krone1999} by Krone. In \cite{Krone1999}, a duality relationship between the two-stage contact process and a `on-off' process is given. Several important open questions are proposed at the end of \cite{Krone1999}, some of which are answered in \cite{Foxall2015}. For instance, it is shown in \cite{Foxall2015} that the complete convergence theorem holds for the two-stage contact process, i.e., the process converges weakly to the convex combination of two invariant distributions.

When $\gamma=+\infty$, i.e., a semi-infected vertex becomes a fully-infected one immediately, the two-stage contact process reduces to the basic contact process introduced in \cite{Har1974}. For a detailed survey about the study of the basic contact process, see Chapter six of \cite{Lig1985} and Part one of \cite{Lig1999}.

\section{Main results}\label{section two}
In this section we give our main results. First we introduce some notations and definitions. For any $t\geq 0$, we define
\[
C_t=\big\{x\in \mathbb{Z}^d:~\eta_t(x)=2\big\}
\]
as the set of fully-infected vertices at the moment $t$ and
\[
D_t=\big\{x\in \mathbb{Z}^d:~\eta_t(x)=1\big\}
\]
as the set of semi-infected vertices at the moment $t$ while $I_t=C_t\bigcup D_t$ as the set of infected vertices at the moment $t$. For $C,D\subseteq \mathbb{Z}^d$, we write $\eta_t,C_t,D_t,I_t$ as $\eta_t^{(C,D)},C_t^{(C,D)},D_t^{(C,D)}, I_t^{(C,D)}$
when $C_0=C, D_0=D$. If $C=\{x\}$ (resp. $D=\{x\}$) for some $x\in \mathbb{Z}^d$, we write $(C,D)$ as $(x,D)$ (resp. $(C,x)$) instead of $(\{x\},D)$ (resp. $(C,\{x\})$). Throughout this paper, we assume that $\delta, \gamma$ are fixed positive constants. We use $P_\lambda$ to denote the probability measure of the two-stage contact process with infection rate $\lambda$. The expectation with respect to $P_\lambda$ is denoted by $E_\lambda$. We write $P_\lambda, E_\lambda$ as $P_{\lambda,d}, E_{\lambda,d}$ when we need to point out the dimension $d$ of the lattice.

It is obviously that $P_\lambda\big(I_t^{(O,\emptyset)}\neq \emptyset \text{~for all~}t\geq 0\big)$ is increasing with $\lambda$, then it is reasonable to define
\begin{equation}\label{equ 2.1 first definition of CV}
\lambda_c=\sup\big\{\lambda:~P_\lambda\big(I_t^{(O,\emptyset)}\neq \emptyset \text{~for all~}t\geq 0\big)=0\big\}.
\end{equation}
$\lambda_c$ is called the critical value of the infection rate. When $\lambda<\lambda_c$, the infected vertices of the two-stage contact process with infection rate $\lambda$ die out with probability one conditioned on $O$ is the unique initially fully-infected vertex while other vertices are healthy at $t=0$.

It is shown in \cite{Krone1999} that the two-stage contact process $\{\eta_t\}_{t\geq 0}$ is a monotonic process with respect to the partial order $\preceq$ on $\{0,1,2\}^{\mathbb{Z}^d}$ that $\eta\preceq\xi$ when and only when $\eta(x)\leq \xi(x)$ for all $x\in \mathbb{Z}^d$. As a result, $\eta_t^{(\mathbb{Z}^d,\emptyset)}$ converges weakly to an invariant distribution $\nu$ as $t\rightarrow+\infty$. $\nu$ is called the upper invariant distribution of the two-stage contact process. We write $\nu$ as $\nu_\lambda$ when we need to point out the infection rate $\lambda$ and further write $\nu_\lambda$ as $\nu_{\lambda,d}$ when we need to point out the dimension $d$ of the lattice.

It is obviously that $\nu_\lambda(\eta(O)\neq 0)$ is increasing with $\lambda$, so it is reasonable to define
\begin{equation}\label{equ 2.2 second definition of CV}
\widetilde{\lambda}_c=\sup\big\{\lambda:~\nu_\lambda(\eta(O)\neq 0)=0\big\}.
\end{equation}

The following proposition is shown in \cite{Foxall2015}.
\begin{proposition}\label{proposition 2015}
(Foxall, 2015) If $\lambda_c$ and $\widetilde{\lambda}_c$ are defined as in Equations \eqref{equ 2.1 first definition of CV} and \eqref{equ 2.2 second definition of CV} respectively, then
\[
\lambda_c=\widetilde{\lambda}_c.
\]
\end{proposition}
Proposition \ref{proposition 2015} shows that the above two types of critical values of the two-stage contact process are equal. So from now on, we use $\lambda_c$ to denote both of them.

We write $\lambda_c$ as $\lambda_c(d)$ when we need to point out the dimension $d$ of the lattice $\mathbb{Z}^d$. It is shown in \cite{Xue2017} that
\[
\lim_{d\rightarrow+\infty}2d\lambda_c(d)=\frac{1+\delta+\gamma}{\gamma}.
\]
As a result, for sufficiently large $d$ and $\lambda>\frac{1+\delta+\gamma}{\gamma}$,
\begin{equation*}
\nu_{\frac{\lambda}{2d},d}(\eta(O)\neq 0)>0.
\end{equation*}
One of our main results in this paper gives a more precise result than the above inequality. To give this result, we define
\[
\pi(A,B)=\nu\Big(\eta(x)\neq 2\text{~for any~}x\in A \text{~and~}\eta(y)=0\text{~for any~}y\in B\Big)
\]
for any $A,B\subseteq \mathbb{Z}^d$ that $A\bigcap B=\emptyset$. We write $\pi(A,B)$ as $\pi(A,B,\lambda,d)$ when we need to point out the infection rate $\lambda$ and the dimension $d$ of the lattice. Then, for any $d\geq 1, m,n\geq 0$ and $\lambda>\frac{1+\delta+\gamma}{\gamma}$, we define
\begin{align*}
&\Pi(m,n,\lambda,d)=\\
&\sup\Bigg\{\bigg|\pi(A,B,\frac{\lambda}{2d},d)-\Big(1-\frac{\lambda\gamma-(1+\delta+\gamma)}{\lambda(\gamma+1)}\Big)^m\Big(\frac{1+\delta+\gamma}
{\lambda\gamma}\Big)^n\bigg|:\\
&~A,B\subseteq\mathbb{Z}^d, |A|=m, |B|=n, A\bigcap B=\emptyset\Bigg\},
\end{align*}
where $|A|$ is the cardinality of $A$. Then, we obtain the following theorem, which is our first main result.
\begin{theorem}\label{theorem 2.1 new main}
For any $\lambda>\frac{1+\delta+\gamma}{\gamma}$ and integers $m,n\geq 0$,
\[
\lim_{d\rightarrow+\infty}\Pi(m,n,\lambda,d)=0.
\]
\end{theorem}
Intuitively, Theorem \ref{theorem 2.1 new main} shows that $\nu_{\frac{\lambda}{2d},d}$ with $\lambda>\frac{1+\delta+\gamma}{\gamma}$ and large $d$ is approximate to the a product measure $m$ on $\{0,1,2\}^{\mathbb{Z}^d}$ that $\{\eta(x):~x\in \mathbb{Z}^d\}$ are independent under $m$ and
\begin{align*}
&m(\eta(x)=0)=\frac{1+\delta+\gamma}{\lambda\gamma},\text{~}m(\eta(x)=2)=\frac{\lambda\gamma-(1+\delta+\gamma)}{\lambda(\gamma+1)}\\
&\text{while~}m(\eta(x)=1)=\frac{\lambda\gamma-(1+\delta+\gamma)}{\lambda\gamma(\gamma+1)}
\end{align*}
for each $x\in \mathbb{Z}^d$.

When $\gamma=+\infty$, the process reduces to the basic contact process. Let $\widetilde{\nu}$ be the upper invariant measure of the basic contact process, then it is shown in \cite{Schonmann1986} that
\begin{align}\label{equ 2.3 plus Schonmann}
\lim_{d\rightarrow+\infty}\sup\Bigg\{\bigg|\widetilde{\nu}_{\frac{\lambda}{2d},d}\Big(\eta(x)=0\text{~for all~}x\in A\Big)
&-\big(\frac{1}{\lambda}\big)^m\bigg|:\notag\\
&~A\subseteq\mathbb{Z}^d,~|A|=m\Bigg\}=0
\end{align}
for each $m\geq 0$ and $\lambda>1$. Since $\lim_{\gamma\rightarrow+\infty}\frac{1+\delta+\gamma}{\lambda\gamma}=\frac{1}{\lambda}$, Theorem \ref{theorem 2.1 new main} can be considered as an extension of Equation \eqref{equ 2.3 plus Schonmann}.

Our second main result is about the asymptotic behavior of $\lambda_c(d)$ as $d\rightarrow+\infty$. According to the approach introduced in \cite{Xue2017},
\[
0\leq 2d\lambda_c(d)-\frac{1+\delta+\gamma}{\gamma}\leq O\Big(\frac{(\log d)^{\frac{3}{\log d}}}{\log d}\Big)
\]
as $d$ grows to infinity. The following theorem gives a stronger conclusion that $2d\lambda_c(d)-\frac{1+\delta+\gamma}{\gamma}$ and $1/d$ are infinitesimals in the same order as $d\rightarrow+\infty$, which is our second main result.

\begin{theorem}\label{theorem 2.1 main}
If $\lambda_c$ is defined as in Equation \eqref{equ 2.1 first definition of CV}, then
\[
f_1\leq \liminf_{d\rightarrow+\infty}d\Big(2d\lambda_c(d)-\frac{1+\delta+\gamma}{\gamma}\Big) \leq \limsup_{d\rightarrow+\infty}d\Big(2d\lambda_c(d)-\frac{1+\delta+\gamma}{\gamma}\Big)\leq f_2,
\]
where
\[
f_1=\frac{1}{2}(1+\frac{1}{\gamma})\frac{(1+\delta+\gamma)^2}{\gamma(2+\delta+\gamma)} \text{\quad and \quad} f_2=\frac{1+\gamma+\delta}{\gamma}(1+\frac{1}{\gamma}),
\]
which are constants only depend on $\gamma$ and $\delta$.
\end{theorem}

The counterpart of Theorem \ref{theorem 2.1 main} for the critical value of the basic contact process is obtained in References \cite{Grif1983, Hol1981} and \cite{Lig1985}. According to the results given in these references, the critical value $\beta_c$ of the basic contact process on $\mathbb{Z}^d$ satisfies
\begin{equation}\label{equ 2.3}
\frac{1}{2}\leq \liminf_{d\rightarrow+\infty}d\Big(2d\beta_c(d)-1\Big)\leq \limsup_{d\rightarrow+\infty}d\Big(2d\beta_c(d)-1\Big)\leq 1.
\end{equation}
The conclusion that $\liminf_{d\rightarrow+\infty}d\Big(2d\beta_c(d)-1\Big)\geq \frac{1}{2}$ follows from a stronger result that $\beta_c(d)\geq \frac{1}{2d-1}$, which is shown in Section 3.5 of \cite{Lig1985}.

Note that $\lim_{\gamma\rightarrow+\infty}f_1(\gamma)=\frac{1}{2}$ while $\lim_{\gamma\rightarrow+\infty}f_2(\gamma)=1$, hence Theorem \ref{theorem 2.1 main} can be considered as an extension of Equation \eqref{equ 2.3}.

It is natural to ask whether there exists $f_3$ such that
\[
\lim_{d\rightarrow+\infty}d\Big(2d\lambda_c(d)-\frac{1+\delta+\gamma}{\gamma}\Big)=f_3.
\]
This question is open even for the basic contact process, i.e, the case where $\gamma=+\infty$. We will work on this question as a further study.

The remainder of this paper is devoted to the proofs of Theorems \ref{theorem 2.1 new main} and \ref{theorem 2.1 main}. Since the proof of Theorem \ref{theorem 2.1 new main} relies on some results occurring in the proof of Theorem \ref{theorem 2.1 main}, we will first prove Theorem \ref{theorem 2.1 main} in
Sections \ref{section three} and \ref{section four}. In Section \ref{section three}, we will prove
\begin{equation}\label{equ 2.4}
\liminf_{d\rightarrow+\infty}d\Big(2d\lambda_c(d)-\frac{1+\delta+\gamma}{\gamma}\Big)\geq f_1.
\end{equation}
The proof of Equation \eqref{equ 2.4} relies on a graphic representation of the two-stage contact process. In Section \ref{section four}, we will prove
\begin{equation}\label{equ 2.5}
\limsup_{d\rightarrow+\infty}d\Big(2d\lambda_c(d)-\frac{1+\delta+\gamma}{\gamma}\Big)\leq f_2.
\end{equation}
The theory of the linear system introduced in Chapter nine of \cite{Lig1985} is crucial for the proof of Equation \eqref{equ 2.5}. A linear system with state space $\big([0,+\infty)\times[0,+\infty)\big)^{\mathbb{Z}^d}$ will be introduced as an auxiliary model.

The proof of Theorem \ref{theorem 2.1 new main} will be given in Sections \ref{section five} and \ref{section six}. In section \ref{section five}, we will introduce a two-type branching process. If there are $m$ semi-infected individuals and $n$ fully-infected individuals for this branching process at $t=0$, then the survival probability of this branching process is an upper bound of $1-\pi(A,B)$ with $|A|=m$ and $|B|=n$. A duality relationship introduced in \cite{Krone1999} between the two-stage contact process and a so-called `on-off' model will be utilized. For details, see Section \ref{section five}.

In Section \ref{section six}, some lower bounds of $1-\pi(A,B)$ will be given. The linear system introduced in Section \ref{section four} and the duality relationship introduced in \cite{Krone1999} are still crucial for us to give these lower bounds. For details, see section \ref{section six}.

\section{Proof of Equation \eqref{equ 2.4}}\label{section three}
In this section we give the proof of Equation \eqref{equ 2.4}. First we introduce a graphic representation of the two-stage contact process. According to this graphic representation, for given $A,B\subseteq \mathbb{Z}^d$ that $A\bigcap B=\emptyset$, the crowd of processes
\[
\Big\{\{I_t^{(C,D)}\}_{t\geq 0}:~C\subseteq A, D\subseteq B\Big\}
\]
can be coupled under a same probability space. We consider the set $\mathbb{Z}^d\times [0,+\infty)$, i.e, there is a time axis $[0,+\infty)$ on each vertex $x\in \mathbb{Z}^d$. For each $x\in \mathbb{Z}^d$, let $\{Y_x(t)\}_{t\geq 0}$ be a Poisson process with rate one, then we put a `$\Delta$' on $(x,s)$ for each event moment $s$ of $Y_x(\cdot)$. For each $x\in \mathbb{Z}^d$, let $\{W_x(t)\}_{t\geq 0}$ be a Poisson process with rate $\delta$, then we put a `$\ast$' on $(x,r)$ for each event moment $r$ of $W_x(\cdot)$. For each $x\in \mathbb{Z}^d$, let $\{V_x(t)\}_{t\geq 0}$ be a Poisson process with rate $\gamma$, then we put a `$\diamond$' on $(x,u)$ for each event moment $u$ of $V_x(\cdot)$. For any $x,y\in \mathbb{Z}^d$ that $x\sim y$, let $\{U_{(x,y)}(t)\}_{t\geq 0}$ be a Poisson process with rate $\lambda$, then we put a `$\rightarrow$' from $(x,v)$ to $(y,v)$ for each event moment $v$ of $U_{(x,y)}(\cdot)$. We assume that all these Poisson processes are independent. Note that we care about the order of $x$ and $y$, hence $U_{(x,y)}\neq U_{(y,x)}$.

Now assuming that $A,B\subseteq \mathbb{Z}^d$ that $A\bigcap B=\emptyset$, then we put a `$\diamond$' on $(x,0)$ for each $x\in A$. For $x\in A\bigcup B, y\in \mathbb{Z}^d$ and $t>0$, we say that there is an infection path from $(x,0)$ to $(y,t)$ when there exist $n\geq 0$, $x=x_0\sim x_1\sim x_2\sim\ldots\sim x_n=y$ and $0=t_0<t_1<t_2<\ldots<t_n<t_{n+1}=t$ such that the following five conditions all hold.

(1) There is an `$\rightarrow$' from $(x_{i-1},t_i)$ to $(x_i,t_i)$ for all $1\leq i\leq n$.

(2) There exists at least one `$\diamond$' on $\{x_i\}\times[t_i,t_{i+1})$ for all $0\leq i\leq n-1$.

(3) There is no `$\Delta$' on $\{x_i\}\times [t_i,t_{i+1})$ for all $0\leq i\leq n$.

(4) For each $0\leq i\leq n-1$, let
\[
m_i=\inf\{s\in [t_i,t_{i+1}):\text{~there is a `$\diamond$' on~}(x_i,s)\},
\]
then there is no `$\ast$' on $\{x_i\}\times [t_i,m_i)$ for all $0\leq i\leq n$.

(5) Let
\[
m_n=\inf\{s\in [t_n,t):\text{~there is a `$\diamond$' on~}(x_i,s)\},
\]
then there is no `$\ast$' on $\{y\}\times [t_n,m_n)$ if $m_n<+\infty$ while there is no `$\ast$' on $\{y\}\times [t_n,t)$ if $m_n=+\infty$.

Note that condition (2) ensures that $m_i<+\infty$ for $0\leq i\leq n-1$ while $m_n$ may equals $\inf\emptyset=+\infty$, so condition (5) contains two cases.

For $C\subseteq A$, $D\subseteq B$ and $t\geq 0$, we define
\begin{align*}
\widehat{I}_t^{(C,D)}=\big\{y\in \mathbb{Z}^d:&\text{~there is an infection path from~}(x,0)\\
&\text{~to~}(y,t)\text{~for some~}x\in C\bigcup D\big\}.
\end{align*}
According to the theory of the graphical method introduced in \cite{Har1978}, it is easy to check that $\{\widehat{I}^{(C,D)}_t\}_{t\geq 0}$ and $\{I_t^{(C,D)}\}_{t\geq 0}$ are identically distributed, where
\[
I_t^{(C,D)}=\big\{x:~\eta_t^{(C,D)}(x)>0\big\}
\]
defined as in Section \ref{section two}. For readers not familiar with the graphical method, we give an intuitive explanation here. An (semi- of fully-) infected vertex $x$ becomes healthy at the event moment of $Y_x(\cdot)$. If $x$ is semi-infected, it also becomes healthy at the event moment of $W_x(\cdot)$ while becomes fully-infected at the event moment of $V_x(\cdot)$. If $x$ is fully infected while the neighbor $y$ of $x$ is healthy, then $y$ is infected by $x$ to become semi-infected when there is an `$\rightarrow$' from $x$ to $y$. As a result, if there is an infection path from $(x,0)$ to $(y,t)$ for $x\in C\bigcup D$, then for each $i\leq n-1$, $x_i$ is (semi- or fully-) infected at $t_i$ and is full-infected at $m_i$ while maintains fully-infected till $t_{i+1}$ to ensure that $x_{i+1}$ is infected at $t_{i+1}$ for all $0\leq i\leq n-1$. Hence, $y=x_n$ is infected at $t_{n}$. If $m_n<+\infty$, then $y$ becomes fully-infected at $m_n$ and maintains fully-infected till $t$. If $m_n=+\infty$, then $y$ maintains semi-infected till $t$. Therefore,
\[
\widehat{I}_t^{(C,D)}\subseteq \big\{y:~\eta_t^{(C,D)}(y)>0\big\}.
\]
The opposite direction that $\widehat{I}_t^{(C,D)}\supseteq \big\{y:~\eta_t^{(C,D)}(x)>0\big\}$ follows from similar analysis, we omit the details.

From now on, we assume that $\Big\{\{I_t^{(C,D)}\}_{t\geq 0}:~C\subseteq A, D\subseteq B\Big\}$ are coupled under a same probability space such that
\begin{align}\label{equ 2.6}
I_t^{(C,D)}=\big\{y\in \mathbb{Z}^d:&\text{~there is an infection path from~}(x,0)\\
&\text{~to~}(y,t)\text{~for some~}x\in C\bigcup D\big\}\notag
\end{align}
for any $t>0$. According to Equation \eqref{equ 2.6}, we have the following lemma, which is crucial for us to prove Equation \eqref{equ 2.4}.
\begin{lemma}\label{lemma 3.1}
For $A,B\subseteq \mathbb{Z}^d$ that $A\bigcap B=\emptyset$ and $C_+,C_-\subseteq A$ while $D_+,D_-\subseteq B$,
\begin{align*}
&P_\lambda\Big(I_t^{(C_+\bigcup C_-,D_+\bigcup D_-)}\neq \emptyset\Big)+P_\lambda\Big(I_t^{(C_+\bigcap C_-,D_+\bigcap D_-)}\neq\emptyset\Big)\\
&\leq P_\lambda\Big(I_t^{(C_+,D_+)}\neq \emptyset\Big)+P_\lambda\Big(I_t^{(C_-,D_-)}\neq\emptyset\Big)
\end{align*}
for any $t\geq 0$.
\end{lemma}

\proof
For $C\subseteq A$ and $D\subseteq B$, we use $H_t(C,D)$ to denote the indicator function of the event that there exists an infection path from $(x,0)$ to $(y,0)$ for some $x\in C\bigcup D$ and $y\in \mathbb{Z}^d$, then its easy to check that
\begin{align*}
&H_t(C_+\bigcup C_-,D_+\bigcup D_-)+H_t(C_+\bigcap C_-,D_+\bigcap D_-)\\
&\leq H_t(C_+,D_+)+H_t(C_-,D_-)
\end{align*}
for $C_+, C_-\subseteq A$ and $D_+,D_-\subseteq B$. Lemma \ref{lemma 3.1} follows from the above inequality directly since
\[
P_\lambda\Big(I_t^{(C,D)}\neq \emptyset\Big)=E_\lambda\Big(H_t(C,D)\Big)
\]
according to Equation \eqref{equ 2.6}.

\qed

For simplicity, we define
\begin{align*}
&\alpha=P_\lambda\Big(I_t^{(\emptyset,O)}\neq \emptyset\text{~for all~}t>0\Big),\\
&q_1=P_\lambda\Big(I_t^{(O,\emptyset)}\neq \emptyset\text{~for all~}t>0\Big),\\
&k_1=P_\lambda\Big(I_t^{(O,e_1)}\neq \emptyset\text{~for all~}t>0\Big),\\
&k_2=\sup\Big\{P_\lambda\Big(I_t^{(O,\{e_1,y\})}\neq \emptyset\text{~for all~}t>0\Big):~y\sim O, y\neq e_1\Big\},\\
&q_2=P_\lambda\Big(I_t^{(\{O,e_1\},\emptyset)}\neq \emptyset\text{~for all~}t>0\Big),\\
&q_3=\sup\Big\{P_\lambda\Big(I_t^{(\{O,e_1\},y)}\neq \emptyset\text{~for all~}t>0\Big):~y\sim O, y\neq e_1\Big\},
\end{align*}
where $e_1=(1,0,\ldots,0)$ defined as in Section \ref{section one},
then we have the following lemma.
\begin{lemma}\label{lemma 3.2}
\[
k_2\leq 2k_1-q_1 \text{\quad and\quad}q_3\leq k_1+q_2-q_1.
\]
\end{lemma}

\proof

For $y\sim O$ and $y\neq \emptyset$, let $C_+=C_-=A=\{O\}$, $B=\{e_1,y\}$, $D_+=\{e_1\}$ and $D_-=\{y\}$, then by Lemma \ref{lemma 3.1} and the spatial homogeneity of the process,
\begin{align*}
&P_\lambda\Big(I_t^{(O,\{e_1,y\})}\neq \emptyset\Big)+P_\lambda\Big(I_t^{(O,\emptyset)}\neq\emptyset\Big)\\
&\leq P_\lambda\Big(I_t^{(O,e_1)}\neq \emptyset\Big)+P_\lambda\Big(I_t^{(O,y)}\neq\emptyset\Big)=2P_\lambda\Big(I_t^{(O,e_1)}\neq \emptyset\Big)
\end{align*}
for any $t\geq 0$. Let $t\rightarrow+\infty$ and then
\[
P_\lambda\Big(I_t^{(O,\{e_1,y\})}\neq \emptyset\text{~for all~}t>0\Big)+q_1\leq 2k_1,
\]
$k_2\leq 2k_1-q_1$ follows from which directly.

Let $C_+=A=\{O,e_1\}$, $C_-=\{O\}$, $B=\{y\}$, $D_+=\emptyset$ and $D_-=\{y\}$, then $q_3\leq k_1+q_2-q_1$ follows from the same analysis, we omit the details.

\qed

Now we give the proof of Equation \eqref{equ 2.4}.

\proof[Proof of Equation \eqref{equ 2.4}]
Let $C_0=\{x:\eta_0(x)=2\}$ and $D_0=\{x:\eta_0(x)=1\}$ defined as in Section \ref{section two}.
If $C_0=\{O\}$ and $D_0=\emptyset$, then according to the property of independent exponential times, a neighbor of $O$ is infected to become semi-infected with probability $\frac{2d\lambda}{1+2d\lambda}$ while $O$ becomes healthy without infecting any neighbor with probability $\frac{1}{1+2d\lambda}$. Therefore, according to the strong Markov property and the spatial homogeneity of the process,
\begin{equation}\label{equ 3.2}
q_1=\frac{2d\lambda}{2d\lambda+1}k_1.
\end{equation}
If $C_0=\emptyset$ and $D_0=\{O\}$, then according to a similar analysis,
\begin{equation}\label{equ 3.3}
\alpha=\frac{\gamma}{1+\delta+\gamma}q_1.
\end{equation}
If $C_0=\{O\}$ and $D_0=\{e_1\}$, then $(C_0,D_0)$ jumps to $(C,D)$ with probability
\[
\begin{cases}
\frac{1+\delta}{(2d-1)\lambda+2+\delta+\gamma} & \text{~if~} C=\{O\}, D=\emptyset,\\
\frac{1}{(2d-1)\lambda+2+\delta+\gamma} & \text{~if~} C=\emptyset, D=\{e_1\},\\
\frac{\gamma}{(2d-1)\lambda+2+\delta+\gamma} & \text{~if~} C=\{O,e_1\}, D=\emptyset,\\
\frac{\lambda}{(2d-1)\lambda+2+\delta+\gamma} & \text{~if~} y\sim O, y\neq e_1, C=\{O\}, D=\{e_1,y\},\\
 0 & \text{~else}.
\end{cases}
\]
Then, according to the strong Markov property, spatial homogeneity of the process and Lemma \ref{lemma 3.2},
\begin{align}\label{equ 3.4}
k_1\leq &\frac{1+\delta}{(2d-1)\lambda+2+\delta+\gamma}q_1+\frac{1}{(2d-1)\lambda+2+\delta+\gamma}\alpha\\
&+\frac{\gamma}{(2d-1)\lambda+2+\delta+\gamma}q_2+\frac{(2d-1)\lambda}{(2d-1)\lambda+2+\delta+\gamma}(2k_1-q_1). \notag
\end{align}
If $C_0=\{O,e_1\}$ and $D_0=\emptyset$, then according to Lemma \ref{lemma 3.2} and a similar analysis with that leads to Equation \eqref{equ 3.4},
\begin{equation*}
q_2\leq \frac{2}{2(2d-1)\lambda+2}q_1+\frac{2(2d-1)\lambda}{2(2d-1)\lambda+2}(k_1+q_2-q_1)
\end{equation*}
and hence
\begin{equation}\label{equ 3.5}
q_2\leq \frac{(4d-1)\lambda}{2d\lambda+1}k_1.
\end{equation}
By Equations \eqref{equ 3.2}, \eqref{equ 3.3}, \eqref{equ 3.4} and \eqref{equ 3.5},
\begin{equation}\label{equ 3.6}
k_1\big(M(\lambda,\delta,\gamma)-1\big)\geq 0,
\end{equation}
where
\begin{align*}
&M(\lambda,\delta,\gamma)=\\
&\frac{\gamma}{(2d-1)\lambda+2+\delta+\gamma}\frac{(4d-1)\lambda}{2d\lambda+1}
+\frac{1}{(2d-1)\lambda+2+\delta+\gamma}\frac{\gamma}{1+\delta+\gamma}\frac{2d\lambda}{2d\lambda+1}\\
&+\frac{1+\delta}{(2d-1)\lambda+2+\delta+\gamma}\frac{2d\lambda}{2d\lambda+1}+\frac{(2d-1)\lambda}{(2d-1)\lambda+2+\delta+\gamma}\frac{2d\lambda+2}{2d\lambda+1}.
\end{align*}
By direct calculation, it is easy to check that
\[
M(\lambda,\delta,\gamma)<1
\]
when
\[
\lambda<\frac{1+\delta+\gamma}{2d\gamma}\frac{2+\delta+\gamma}{1+[1-\frac{1}{2d}(1+\frac{1}{\gamma})](1+\delta+\gamma)}.
\]
Note that here we assume that $d$ is sufficiently large that
\[
1+[1-\frac{1}{2d}(1+\frac{1}{\gamma})](1+\delta+\gamma)>0.
\]
Then, by Equations \eqref{equ 3.6} and \eqref{equ 3.2}, $k_1=0$ and hence $q_1=0$
when
\[
\lambda<\frac{1+\delta+\gamma}{2d\gamma}\frac{2+\delta+\gamma}{1+[1-\frac{1}{2d}(1+\frac{1}{\gamma})](1+\delta+\gamma)}.
\]
Then, according to the definition of $\lambda_c(d)$ given in Equation \eqref{equ 2.1 first definition of CV},
\begin{equation}\label{equ 3.7}
\lambda_c(d)\geq \frac{1+\delta+\gamma}{2d\gamma}\frac{2+\delta+\gamma}{1+[1-\frac{1}{2d}(1+\frac{1}{\gamma})](1+\delta+\gamma)}.
\end{equation}
Since
\[
\frac{2+\delta+\gamma}{1+[1-\frac{1}{2d}(1+\frac{1}{\gamma})](1+\delta+\gamma)}
=1+\frac{1}{2d}(1+\frac{1}{\gamma})\frac{1+\delta+\gamma}{2+\delta+\gamma}+O(\frac{1}{d^2}),
\]
Equation \eqref{equ 2.4} follows from Equation \eqref{equ 3.7} directly.

\qed

\section{Proof of Equation \eqref{equ 2.5}}\label{section four}
In this section we give the proof of Equation \eqref{equ 2.5}. First we introduce a continuous-time
Markov process $\{\rho_t\}_{t\geq 0}$ as an auxiliary process for the proof.
The state space of $\{\rho_t\}_{t\geq 0}$ is $X_2=\big([0,+\infty)\times [0,+\infty)\big)^{\mathbb{Z}^d}$, i.e., at each vertex
$x\in \mathbb{Z}^d$ there is a spin $\rho(x)=\big(\zeta(x),g(x)\big)$ that $\zeta(x),g(x)\geq 0$. The generator function of $\{\rho_t\}_{t\geq 0}$ is given by
\begin{align}\label{equ 4.1 generator of linear system}
\Omega_2 f(\rho)=&\sum_{x\in \mathbb{Z}^d}\big[f(\rho^x)-f(\rho)\big]+\delta\sum_{x\in \mathbb{Z}^d}\big[f(\rho^{x,+})-f(\rho)\big] \\
&+\gamma\sum_{x\in \mathbb{Z}^d}\big[f(\rho^{x,-})-f(\rho)\big]+\lambda\sum_{x\in \mathbb{Z}^d}\sum_{y\sim x}\big[f(\rho^{x,y})-f(\rho)\big]\notag
\end{align}
for any $\rho\in X_2$ and $f\in C(X_2)$, where
\[
\rho^x(y)=
\begin{cases}
\rho(y)=\big(\zeta(y),g(y)\big) & \text{~if~}y\neq x,\\
\big(0,0\big) & \text{~if~}y=x,
\end{cases}
\]
\[
\rho^{x,+}(y)=
\begin{cases}
\rho(y)=\big(\zeta(y),g(y)\big) & \text{~if~}y\neq x,\\
\big(\zeta(x),0\big) & \text{~if~}y=x,
\end{cases}
\]
\[
\rho^{x,-}(y)=
\begin{cases}
\rho(y)=\big(\zeta(y),g(y)\big) & \text{~if~}y\neq x,\\
\big(\zeta(x)+\frac{1}{\gamma}g(x),0\big) & \text{~if~}y=x,
\end{cases}
\]
and
\[
\rho^{x,y}(z)=
\begin{cases}
\rho(z)=\big(\zeta(z),g(z)\big) & \text{~if~}z\neq x,\\
\big(\zeta(x),g(x)+b\zeta(y)\big) & \text{~if~}z=x,
\end{cases}
\]
where $b=\frac{1+\delta+\gamma}{2d\lambda}$.

According to the generator $\Omega_2$ of $\{\rho_t\}_{t\geq 0}$, if the state of $\{\rho_t\}_{t\geq 0}$
jumps at a moment $s$, then $\zeta_s(x), g_s(x)$ are linear combinations of $\{\zeta_{s-}(y)\}_{y\in \mathbb{Z}^d}$
and $\{g_{s-}(y)\}_{y\in \mathbb{Z}^d}$ for each $x\in \mathbb{Z}^d$. As a result, $\{\rho_t\}_{t\geq 0}$ can be considered as a
linear system, the theory of which is introduced in Chapter nine of \cite{Lig1985}.

In Chapter nine of \cite{Lig1985}, the state space of a linear
system is defined to be $[0,+\infty)^S$, where $S$ is a countable
set. Note that $\{\rho_t\}_{t\geq 0}$ is consistent with this
definition since
$\big([0,+\infty)\times[0,+\infty)\big)^{\mathbb{Z}^d}$ can be
identified with $[0,+\infty)^{\mathbb{Z}^d\times\{1,2\}}$ while
$\mathbb{Z}^d\times \{1,2\}$ is a countable set.

The intuitive explanation of the evolution of $\{\rho_t\}_{t\geq 0}$ is as follows. For any $x\in \mathbb{Z}^d$, its state $\rho(x)=\big(\zeta(x),g(x)\big)$ flips to $(0,0)$ at rate $1$, to $\big(\zeta(x),0\big)$ at rate $\delta$, to $\big(\zeta(x)+\frac{1}{\gamma}g(x),0\big)$ at rate $\gamma$ or to $\big(\zeta(x),g(x)+b\zeta(y)\big)$ at rate $\lambda$ for each neighbor $y$.

From now on, we assume that $\rho_0(x)=(1,1)$ for any $x\in
\mathbb{Z}^d$, then we define
\[
\phi_t(x)=
\begin{cases}
2 & \text{~if~} \zeta_t(x)>0,\\
1 & \text{~if~} \zeta_t(x)=0 \text{~and~}g_t(x)>0,\\
0 & \text{~if~} \zeta_t(x)=g_t(x)=0
\end{cases}
\]
for any $t\geq 0$ and $x\in \mathbb{Z}^d$, where
$\rho_t(x)=\big(\zeta_t(x),g_t(x)\big)$. The following lemma gives the relationship between $\{\rho_t\}_{t\geq 0}$ and the two-stage contact process.
\begin{lemma}\label{lemma 4.1}
$\{\phi_t\}_{t\geq 0}$ is a version of the two-stage contact process with generator \eqref{equ 1.1 generator}.
\end{lemma}

\proof

We only need to check that $\{\phi_t\}_{t\geq 0}$ evolves in the same way as that of the two-stage contact process $\{\eta_t\}_{t\geq 0}$. For any $x\in \mathbb{Z}^d$, if $\phi(x)=0$, i.e., $\rho(x)=\big(0,0\big)$, then $\phi(x)$ flips to $1$ when and only when $\rho(x)$ flips to
\[
\big(0,0+b\zeta(y)\big)=\big(0,b\zeta(y)\big)
\]
for some $y\sim x$ that $\zeta(y)>0$, i.e, $\phi(y)=2$. Since $\rho(x)$ flips to $\big(\zeta(x),g(x)+b\zeta(y)\big)$ at rate $\lambda$, $\phi(x)$ flips from $0$ to $1$ at rate
\[
\lambda\sum_{y\sim x}1_{\{\phi(y)=2\}}=H(x,1,\phi)
\]
defined in Equation \eqref{equ 1.1 generator}. Through a similar way, it is easy to check that in every case $\phi(x)$ flips to a different state $i$ at rate $H(x,i,\phi)$ defined in Equation \eqref{equ 1.1 generator} and the proof is complete.

\qed

By Lemma \ref{lemma 4.1}, from now on we assume that $\{\rho_t\}_{t\geq 0}$ and the two-stage contact process $\{\eta_t\}_{t\geq 0}$ are defined under the same probability space such that
\[
\eta_t(x)=
\begin{cases}
2 & \text{~if~} \zeta_t(x)>0,\\
1 & \text{~if~} \zeta_t(x)=0 \text{~and~}g_t(x)>0,\\
0 & \text{~if~} \zeta_t(x)=g_t(x)=0
\end{cases}
\]
for any $t\geq 0$ and $x\in \mathbb{Z}^d$. As a result,
\begin{equation}\label{equ 4.2}
\nu_\lambda(\eta(O)=2)=\lim_{t\rightarrow+\infty}P_\lambda\Big(\eta_t^{(\mathbb{Z}^d,\emptyset)}(O)=2\Big)=
\lim_{t\rightarrow+\infty}P_\lambda\Big(\zeta_t(O)>0\Big).
\end{equation}

By Equation \eqref{equ 4.2}, we have the following lemma about the
upper bound of the critical value $\lambda_c$.

\begin{lemma}\label{lemma 4.2}
If $\lambda$ satisfies
\[
\sup_{t\geq 0}E_\lambda\big(\zeta_t^2(O)\big)<+\infty,
\]
then $\lambda_c\leq \lambda$.
\end{lemma}

\proof By Equation \eqref{equ 4.2} and Cauchy-Schwartz inequality,
\begin{align}\label{equ 4.4}
\lim_{t\rightarrow+\infty}P_\lambda\Big(\eta_t^{(\mathbb{Z}^d,\emptyset)}(O)=2\Big)&=
\lim_{t\rightarrow+\infty}P_\lambda\Big(\zeta_t(O)>0\Big) \notag\\
&\geq \limsup_{t\rightarrow+\infty}
\frac{\big(E_\lambda\zeta_t(O)\big)^2}{E_\lambda\big(\zeta_t^2(O)\big)}.
\end{align}
Let $\{S(t)\}_{t\geq 0}$ be the semi-group of $\{\rho_t\}_{t\geq
0}$. According to Theorem 9.1.27 of \cite{Lig1985}, which is an
extension version of the Hille-Yosida Theorem for the linear system,
we can execute the calculation
\[
\frac{d}{dt}S(t)f(\rho)=S(t)\Omega_2f(\rho)
\]
for $f$ with the form $f(\rho)=\zeta(x)$ and $f(\rho)=g(x)$. Then,
according to the definition of $\Omega_2$,
\begin{equation*}
\begin{cases}
&\frac{d}{dt}E_\lambda\zeta_t(O)=-E_\lambda\zeta_t(O)+\gamma
E_\lambda\big(\frac{1}{\gamma}g_t(O)\big),\\
&\frac{d}{dt}E_\lambda g_t(O)=-(1+\delta+\gamma)E_\lambda
g_t(O)+\lambda b\sum_{y:y\sim O}E_\lambda\zeta_t(y)
\end{cases}
\end{equation*}
and hence
\begin{equation}\label{equ 4.5}
\begin{cases}
&\frac{d}{dt}E_\lambda\zeta_t(O)=-E_\lambda\zeta_t(O)+
E_\lambda g_t(O),\\
&\frac{d}{dt}E_\lambda g_t(O)=-(1+\delta+\gamma)E_\lambda
g_t(O)+(1+\delta+\gamma)E_\lambda\zeta_t(O)
\end{cases}
\end{equation}
by the spatial homogeneity of $\{\rho_t\}_{t\geq 0}$.

Since $\zeta_0(O)=g_0(O)=1$, it is easy to check that the unique
solution to ODE \eqref{equ 4.5} is
\[
E_\lambda \zeta_t(O)=E_\lambda g_t(O)\equiv 1.
\]
Then, by Equation \eqref{equ 4.4},
\begin{align}\label{equ 4.6}
\nu_\lambda\big(\eta(O)=2\big)&=\lim_{t\rightarrow+\infty}P_\lambda\Big(\eta_t^{(\mathbb{Z}^d,\emptyset)}(O)=2\Big)\\
&\geq
\limsup_{t\rightarrow+\infty}
\frac{1}{E_\lambda\big(\zeta_t^2(O)\big)}\geq \frac{1}{\sup_{t\geq 0}E_\lambda\big(\zeta_t^2(O)\big)}>0
\notag
\end{align}
if $\lambda$ satisfies
\[
\sup_{t\geq 0}E_\lambda\big(\zeta_t^2(O)\big)<+\infty.
\]
Lemma \ref{lemma 4.2} follows directly from Equation \eqref{equ 4.6}
and the equivalent definition of $\lambda_c$ given in Equation
\eqref{equ 2.2 second definition of CV}.

\qed

By Lemma \ref{lemma 4.2}, we want to bound $E_\lambda\big(\zeta_t^2(O)\big)$ from above. For this purpose, we define
\[
F_t(x,1)=E_\lambda\big(\zeta_t(O)\zeta_t(x)\big), F_t(x,2)=E_\lambda\big(\zeta_t(O)g_t(x)\big)
\]
and $F_t(x,3)=E_\lambda\big(g_t(O)g_t(x)\big)$ for each $x\in \mathbb{Z}^d$ and any $t\geq 0$. For any $t>0$, we define
\[
F_t=\Big\{F_t(x,i):~x\in \mathbb{Z}^d, i\in\{1,2,3\}\Big\}
\]
as a function on $X_4=\mathbb{Z}^d\times \{1,2,3\}$. For a $X_4\times X_4$ matrix
\[
G=\{G\big((x,i),(y,j)\big)\}_{(x,i),(y,j)\in X_4}
\]
and two functions $F_+, F_-$ on $X_4$, we write
\[
F_+=GF_-
\]
when and only when
\[
F_+(x,i)=\sum_{(y,j)\in X_4}G\big((x,i),(y,j)\big)F_-(y,j)
\]
for any $(x,i)\in X_4$, as the product of finite dimensional matrixes. Then we have the following lemma.
\begin{lemma}\label{lemma 4.3}
Let
\[
\frac{d}{dt}F_t=\Big\{\frac{d}{dt}F_t(x,i):~(x,i)\in X_4\Big\},
\]
then
\begin{equation}\label{equ 4.7 two}
\frac{d}{dt}F_t=G_\lambda F_t,
\end{equation}
where $G_\lambda$ is a $X_4\times X_4$ matrix such that
\[
G_\lambda\big((x,i),(y,j)\big)=
\begin{cases}
-2 & \text{~if~} x\neq O, i=1 \text{~and~}(y,j)=(x,1),\\
2 & \text{~if~} x\neq O, i=1 \text{~and~}(y,j)=(x,2),\\
-1 & \text{~if~} x=O, i=1 \text{~and~}(y,j)=(O,1),\\
2 & \text{~if~} x=O, i=1 \text{~and~}(y,j)=(O,2),\\
\frac{1}{\gamma} & \text{~if~} x=O,i=1 \text{~and~}(y,j)=(O,3),\\
-(2+\delta+\gamma) & \text{~if~}x\neq O, i=2 \text{~and~}(y,j)=(x,2),\\
1 & \text{~if~}x\neq O, i=2 \text{~and~}(y,j)=(x,3),\\
\frac{1+\delta+\gamma}{2d} & \text{~if~}x\neq O, i=2, y\sim x \text{~and~}j=1,\\
-(1+\delta+\gamma) & \text{~if~}x=O, i=2 \text{~and~}(y,j)=(O,2),\\
1+\delta+\gamma & \text{~if~}x=O, i=2 \text{~and~}(y,j)=(e_1,1),\\
-2(1+\delta+\gamma) & \text{~if~}x\neq O, i=3 \text{~and~}(y,j)=(x,3),\\
\frac{1+\delta+\gamma}{d} & \text{~if~}x\neq O, i=3, y\sim x\text{~and~}j=2,\\
-(1+\delta+\gamma) & \text{~if~}x=O, i=3 \text{~and~}(y,j)=(O,3),\\
2(1+\delta+\gamma) & \text{~if~}x=O, i=3 \text{~and~}(y,j)=(e_1,2),\\
\frac{(1+\delta+\gamma)^2}{2d\lambda} & \text{~if~}x=O,i=3 \text{~and~} (y,j)=(O,1),\\
0 & \text{~else}.
\end{cases}
\]
\end{lemma}

\proof

According to the spatial homogeneity of the process
$\{\rho_t\}_{t\geq 0}$,
\begin{align}\label{equ 4.7}
&E_\lambda\big(\zeta_t(u)\zeta_t(v)\big)=F_t(u-v,1)=F_t(v-u,1),\notag\\
&E_\lambda\big(\zeta_t(u)g_t(v)\big)=E_\lambda\big(\zeta_t(v)g_t(u)\big)=F_t(u-v,2)=F_t(v-u,2),\\
&E_\lambda\big(g_t(u)g_t(v)\big)=F_t(u-v,3)=F_t(v-u,3),\notag\\
&F_t(e_1,i)=F_t(y,i)\notag
\end{align}
for any $y\sim O$, $u,v\in \mathbb{Z}^d$ and $i\in \{1,2,3\}$.
Theorem 9.3.1 of \cite{Lig1985} is an extension version of the
Hille-Yosida Theorem for the linear system, according to which we
can execute the calculation that
\[
\frac{d}{dt}S(t)f(\rho)=S(t)\Omega_2f(\rho)
\]
for $f$ with form $f(\rho)=\zeta(x)\zeta(y)$, $f(\rho)=\zeta(x)g(y)$
and $f(\rho)=g(x)g(y)$ for $x,y\in \mathbb{Z}^d$. Therefore, by
Equation \eqref{equ 4.7} and the definition of $\Omega_2$,
\begin{align}\label{equ 4.8}
& \frac{d}{dt}F_t(x,1)=-2F_t(x,1)+2F_t(x,2), \\
& \frac{d}{dt}F_t(x,2)=-(2+\delta+\gamma)F_t(x,2)+F_t(x,3)+\frac{1+\delta+\gamma}{2d}\sum_{y:y\sim x}F_t(y,1), \notag\\
&
\frac{d}{dt}F_t(x,3)=-2(1+\delta+\gamma)F_t(x,3)+\frac{(1+\delta+\gamma)}{d}\sum_{y:y\sim
x}F_t(y,2) \notag
\end{align}
when $x\neq O$ and
\begin{align}\label{equ 4.9}
\frac{d}{dt}F_t(O,1)=&-F_t(O,1)+2F_t(O,2)+\frac{1}{\gamma}F_t(O,3),
\\
\frac{d}{dt}F_t(O,2)=&-(1+\delta+\gamma)F_t(O,2)+(1+\delta+\gamma)F_t(e_1,1),
\notag\\
\frac{d}{dt}F_t(O,3)=&-(1+\delta+\gamma)F_t(O,3)\notag\\
&+2(1+\delta+\gamma)F_t(e_1,2)+\frac{(1+\delta+\gamma)^2}{2d\lambda}F_t(O,1). \notag
\end{align}
Lemma \ref{lemma 4.3} follows from Equations \eqref{equ 4.8} and
\eqref{equ 4.9} directly.

\qed

According to Lemma \ref{lemma 4.3}, we have the following lemma
about a sufficient condition for $\sup_{t\geq 0}
E_\lambda\big(\zeta_t^2(O)\big)<+\infty$.

\begin{lemma}\label{lemma 4.4}
If $\lambda$ satisfies that there exists $K_\lambda:
X_4\rightarrow[0,+\infty)$ such that
\[
G_\lambda K_\lambda=0 \text{~(zero function on $X_4$)}
\]
and $\inf_{(x,i)\in X_4}K_\lambda(x,i)>0$, then
\[
\sup_{t\geq 0} E_\lambda\big(\zeta_t^2(O)\big)<+\infty.
\]
\end{lemma}
To prove Lemma \ref{lemma 4.4}, we need to define the product of two
$X_4\times X_4$ matrixes. For two $X_4\times X_4$ matrixes $G_+$ and
$G_-$, $G_+G_-$ is defined as a $X_4\times X_4$ matrixes that
\[
(G_+G_-)\big((x,i),(y,j)\big)=\sum_{(u,l)\in
X_4}G_+\big((x,i),(u,l)\big)G_-\big((u,l),(y,j)\big)
\]
for any $(x,i), (y,j)\in X_4$, conditioned on the sum is absolute
convergence (otherwise $G_+G_-$ does not exists). Note that this
definition is the same as that of the product of two finite
dimensional matrix, except that the sum must convergence since there
are infinite many terms. Then, we use $G_\lambda^2$ to denote
$G_\lambda G_\lambda$ and define
$G_\lambda^{n+1}=G_\lambda^nG_\lambda$ for $n\geq 2$ by induction.
It is easy to check that the definition of $G_\lambda^n$ is
reasonable for each $n\geq 2$ since for each $(x,i)\in X_4$,
\[
G_\lambda\big((x,i),(y,j)\big)\neq 0
\]
holds for only finite many $(y,j)$s. It is also easy to check that
\[
\sum_{n=0}^{+\infty}\frac{t^n|G_\lambda^n\big((x,i),(y,j)\big)|}{n!}<+\infty
\]
for any $t\geq 0$ and $(x,i), (y,j)\in X_4$, where
$G_\lambda^1=G_\lambda$ and $G_\lambda^0$ is the identity matrix,
then it is reasonable to define $e^{tG_\lambda}$ as the $X_4\times
X_4$ matrix that
\[
e^{tG_\lambda}\big((x,i),(y,j)\big)=\sum_{n=0}^{+\infty}\frac{t^nG_\lambda^n\big((x,i),(y,j)\big)}{n!}
\]
for any $(x,i),(y,j)\in X_4$. Now we can give the proof of Lemma
\ref{lemma 4.4}.

\proof[Proof of Lemma \ref{lemma 4.4}] Since $G_\lambda
K_\lambda=0$, $K_\lambda$ can be considered as the eigenvector of
$G_\lambda$ with respect to the eigenvalue $0$, then according to a
similar analysis with that in the theory of finite-dimensional
linear algebra, $K_\lambda$ is the eigenvector of $e^{tG_\lambda}$
with respect to the eigenvalue $e^{t\times0}=1$, i.e.,
\begin{equation}\label{equ 4.10}
K_\lambda(x,i)=\sum_{(y,j)\in
X_4}e^{tG_\lambda}\big((x,i),(y,j)\big)K_\lambda(y,j)
\end{equation}
for any $t\geq 0$ and $(x,i)\in X_4$.

For any function $K$ on $X_4$, we define
\[
\|K\|_\infty=\sup\big\{|K(x,i)|:~(x,i)\in X_4\big\}
\]
as the $l_\infty$ norm of $K$. Then, we define $X_5$ as the set of
functions on $X_4$ with finite $l_\infty$ norm $\|\cdot\|_{\infty}$.
It is easy to check that $X_5$ is a Banach space with norm
$\|\cdot\|_\infty$. It is also easy to check that there exists a
constant $Q(\lambda)>0$ such that
\[
\|G_\lambda K_+-G_\lambda K_-\|_\infty\leq
Q(\lambda)\|K_+-K_-\|_\infty
\]
for any $K_-,K_+\in X_5$, i.e., ODE \eqref{equ 4.7 two} satisfies
the Lipschitz condition. Since $X_5$ is a Banach space and ODE
\eqref{equ 4.7 two} satisfies the Lipschitz condition, it is easy to
extend the theory of the finite-dimensional linear ODE to the
infinite-dimensional linear ODE \eqref{equ 4.7 two} that the unique
solution to ODE \eqref{equ 4.7 two} is
\[
F_t=e^{tG_\lambda}F_0
\]
i.e.,
\begin{equation}\label{equ 4.12}
F_t(x,i)=\sum_{(y,j)\in
X_4}e^{tG_\lambda}\big((x,i),(y,j)\big)F_0(y,j)
\end{equation}
for any $t\geq 0$ and $(x,i)\in X_4$. Since
$G_\lambda\big((x,i),(y,j)\big)\geq 0$ when $(x,i)\neq (y,j)$, it is
easy to check that
\[
e^{tG_\lambda}\big((x,i),(y,j)\big)\geq 0
\]
for any $(x,i), (y,j)\in X_4$. Then, according to Equations
\eqref{equ 4.10}, \eqref{equ 4.12} and the fact that $F_0(x,i)=1$
for any $(x,i)\in X_4$,
\begin{align}\label{equ 4.13}
E_\lambda\big(\zeta_t(O)\zeta_t(x)\big)&=F_t(x,1)=\sum_{(y,j)\in
X_4}e^{tG_\lambda}\big((x,1),(y,j)\big)  \\
&\leq \sum_{(y,j)\in
X_4}e^{tG_\lambda}\big((x,1),(y,j)\big)\frac{K_\lambda(y,j)}{\inf_{(x,i)\in
X_4} K_\lambda(x,i)} \notag\\
&=\frac{K_\lambda(x,1)}{\inf_{(x,i)\in X_4}K_\lambda(x,i)}<+\infty \notag
\end{align}
for any $t\geq 0$. Let $x=O$, then Lemma \ref{lemma 4.4} follows from Equation
\eqref{equ 4.13} directly.

\qed

By Lemma \ref{lemma 4.4}, we want to find $\lambda$ which ensures the existence of the positive eigenvector $K_\lambda$
of $G_\lambda$ with respect to the eigenvalue $0$. For this purpose, we need two random walks. We denote by $\{S_n\}_{n\geq 0}$ the simple random walk on $\mathbb{Z}^d$ that
\[
P\big(S_{n+1}=y\big|S_n=x\big)=\frac{1}{2d}
\]
for any $n\geq 0$ and $x,y\in \mathbb{Z}^d, x\sim y$. Let $\{\theta_n\}_{n\geq 0}$ be a random walk on
\[
X_4\setminus \{(O,3)\}=\big\{(x,i):~x\in \mathbb{Z}^d, i\in\{1,2,3\}\text{~and}(x,i)\neq (O,3)\big\}
\]
that for each $n\geq 0$,
\begin{align*}
&P\big(\theta_{n+1}=(y,j)\big|\theta_n=(x,i)\big)=\\
&
\begin{cases}
1 & \text{~if~}x\neq O, i=1\text{~and~}(y,j)=(x,2),\\
\frac{1}{2+\delta+\gamma} & \text{~if~}x\neq O, i=2 \text{~and~}(y,j)=(x,3),\\
\frac{1}{2d}\frac{1+\delta+\gamma}{2+\delta+\gamma} & \text{~if~}x\neq O,i=2, y\sim x \text{~and~}j=1,\\
\frac{1}{2d} & \text{~if~}x\neq O, i=3, y\sim x \text{~and~}j=2,\\
1 & \text{~if~}(x,i)=(y,j)=(O,1),\\
1 & \text{~if~}(x,i)=(O,2)  \text{~and~}(y,j)=(e_1,1),\\
0 & \text{~else,}
\end{cases}
\end{align*}
then we define
\[
\Gamma(x,i)=P\Big(\theta_n=(O,1)\text{~for some~}n\geq 0\Big|\theta_0=(x,i)\Big)
\]
for $(x,i)\in X_4$ that $(x,i)\neq (O,3)$, i.e., $\Gamma(x,i)$ is the probability that $\{\theta_n\}_{n\geq 0}$ visits $(O,1)$ at least once conditioned on $\theta_0=(x,i)$. By the definition of $\{\theta_t\}_{t\geq 0}$ and the strong Markov property, $\Gamma(x,i)$ satisfies
\begin{align}\label{equ 4.14}
&\Gamma(x,1)=\Gamma(x,2) \text{~if~} x\neq O,\\
&\Gamma(x,2)=\frac{1}{2+\delta+\gamma+\lambda}\Gamma(x,3)+\frac{1}{2d}\frac{1+\delta+\gamma+\lambda}{2+\delta+\gamma+\lambda}\sum_{y:y\sim x}\Gamma(y,1) \text{~if~}x\neq O, \notag\\
&\Gamma(x,3)=\frac{1}{2d}\sum_{y:y\sim x}\Gamma(y,2) \text{~if~}x\neq O, \notag\\
&\Gamma(O,2)=\Gamma(e_1,1) \text{~and~}\Gamma(O,1)=1.\notag
\end{align}
For any $x\in \mathbb{Z}^d$, we define
\[
\widetilde{\Gamma}(x)=P\Big(S_n=O\text{~for some~}n\geq 0\Big|S_0=x\Big)
\]
as the probability that $\{S_n\}_{n\geq 0}$ visits $O$ at least once conditioned on $S_0=x$. We claim that
\begin{equation}\label{equ 4.15}
\Gamma(x,1)\leq \widetilde{\Gamma}(x)
\end{equation}
for $x\neq O$. Equation \eqref{equ 4.15} follows from the following analysis. For each $n\geq 0$, we write $\theta_n$ as $\big(\theta_n(1),\theta_n(2)\big)$ that $\theta_n(1)\in \mathbb{Z}^d$ and $\theta_n(2)\in \{0,1,2\}$. Conditioned on $\theta_0=(x,1)$ with $x\neq O$, $\{\theta_n(1)\}_{n\geq 0}$ is a lazy version of $\{S_n\}_{n\geq 1}$ with $S_0=x$ until the first moment $n_0$ that $\theta_{n_0}(1)=O$ according to the definition of $\{\theta_n\}_{n\geq 0}$. In other words, before hitting $O$, $\theta(1)$ chooses each neighbor to jump with the same probability $\frac{1}{2d}$ when $\theta(1)$ jumps at some steps while $\theta(1)$ stays still at other steps. Therefore,
\begin{align*}
\Gamma(x,1)&=P\Big(\theta_n=(O,1)\text{~for some~}n\geq 0\Big|\theta_0=(x,1)\Big) \\
&\leq P\Big(\theta_n(1)=O\text{~for some~}n\geq 0\Big|\theta_0=(x,1)\Big)\\
&=P\Big(S_n=O\text{~for some~}n\geq 0\Big|S_0=x\Big)=\widetilde{\Gamma}(x)
\end{align*}
and hence Equation \eqref{equ 4.15} holds. According to the result given in \cite{Kesten1964},
\begin{equation}\label{equ 4.16}
\widetilde{\Gamma}(e_1)=\frac{1}{2d}+\frac{1}{2d^2}+O(\frac{1}{d^3})
\end{equation}
as the dimension $d$ of the lattice grows to infinity.
By Equation \eqref{equ 4.16},
\begin{equation}\label{equ 4.18}
\gamma-(2\gamma+2)\widetilde{\Gamma}(e_1)>0
\end{equation}
when the dimension $d$ of the lattice is sufficiently large. Now we can give the proof of Equation \eqref{equ 2.5}.

\proof[Proof of Equation \eqref{equ 2.5}]

We assume that the dimension $d$ of the lattice is sufficiently large such that Equation \eqref{equ 4.18} holds, then we define
\[
\widetilde{\lambda}=\frac{1+\delta+\gamma}{2d\big[\gamma-(2\gamma+2)\widetilde{\Gamma}(e_1)\big]},
\]
which is positive. Furthermore, we define
\[
h_\lambda=\frac{\gamma[1-2\Gamma(O,2)]-2\Gamma(e_1,2)-\frac{1+\delta+\gamma}{2d\lambda}}{\gamma+2+\frac{1+\delta+\gamma}{2d\lambda}}.
\]
 According to Equation \eqref{equ 4.15} and the fact that $\Gamma(O,2)=\Gamma(e_1,1)$ while
\[
 \Gamma(x,1)=\Gamma(x,2)
\]
for $x\neq O$, it is easy to check that $h_\lambda>0$ when $\lambda>\widetilde{\lambda}$. For any $\lambda>\widetilde{\lambda}$, we define $K_\lambda: X_4\rightarrow[0,+\infty)$ as
\[
K_\lambda(x,i)
=\begin{cases}
\Gamma(x,i)+h_\lambda & \text{~if~} (x,i)\neq (O,3),\\
\gamma[1-2\Gamma(e_1,1)-h_\lambda] & \text{~if~} (x,i)=(O,3).
\end{cases}
\]
Since $h_\lambda\leq 1-2\Gamma(O,2)=1-2\Gamma(e_1,1)$,
\begin{equation}\label{equ 4.19}
\inf_{(x,i)\in X_4}K_\lambda(x,i)\geq \inf\big\{h_\lambda, \gamma[1-2\Gamma(e_1,1)-h_\lambda]\big\}>0
\end{equation}
for $\lambda>\widetilde{\lambda}$. By the definition of $G_\lambda$ and Equation \eqref{equ 4.14}, it is to check that
\[
G_\lambda K_\lambda=0 \text{~(zero function on $X_4$)}
\]
by direct calculation. Then, by Lemma \ref{lemma 4.4} and Equation \eqref{equ 4.19},
\[
\sup_{t\geq 0}E_\lambda\big(\zeta_t^2(O)\big)<+\infty
\]
when $\lambda>\widetilde{\lambda}$. Therefore, by Lemma \ref{lemma 4.2},
\[
\lambda_c\leq \lambda
\]
for any $\lambda>\widetilde{\lambda}$ and hence
\begin{equation}\label{equ 4.20}
\lambda_c\leq \widetilde{\lambda}=\frac{1+\delta+\gamma}{2d\big[\gamma-(2\gamma+2)\widetilde{\Gamma}(e_1)\big]}.
\end{equation}
By utilizing the fact that $\frac{1}{1-x}=\sum_{n=0}^{+\infty}x^n$ for $x\in (0,1)$, Equation \eqref{equ 2.5} follows directly from Equations \eqref{equ 4.16} and \eqref{equ 4.20}.

\qed

\section{Upper bounds of $\text{~}1-\pi(A,B)$}\label{section five}
In this section we will prove the following lemma, which gives upper bounds of $1-\pi(A,B)$.
\begin{lemma}\label{lemma 5.1}
For any $\lambda>\frac{1+\gamma+\delta}{\gamma}$, $d\geq 1$, $m,n\geq 0$ and $A,B\subseteq \mathbb{Z}^d$ that $|A|=m, |B|=n$ while
$A\bigcap B=\emptyset$,
\[
1-\pi(A,B,\frac{\lambda}{2d},d)\leq 1-\bigg(1-\frac{\lambda\gamma-(1+\delta+\gamma)}{\lambda(\gamma+1)}\bigg)^m\bigg(\frac{1+\delta+\gamma}{\lambda\gamma}\bigg)^n.
\]
\end{lemma}
To prove Lemma \ref{lemma 5.1}, we need two auxiliary processes. The first is the `on-off' process introduced in \cite{Krone1999}. The second is a two-type branching process. The `on-off' process $\{\xi_t\}_{t\geq 0}$ is a continuous-time Markov process with state space $\{0,1,2\}^{\mathbb{Z}^d}$ and transition rates function given as follows. For each $x\in \mathbb{Z}^d$ and $t\geq 0$,
\begin{align}\label{equ 5.2}
&\xi_t(x) \text{~flips from $i$ to $j$ at rate~}\\
&\begin{cases}
1 & \text{~if~}i\in \{1,2\}\text{~and~}j=0,\\
\delta & \text{~if~}i=2 \text{~and~}j=1,\\
\gamma & \text{~if~}i=1 \text{~and~}j=2,\\
\lambda\sum_{y\sim x}1_{\{\xi_t(y)=2\}} & \text{~if~}i=0 \text{~and~}j=1,\\
0 & \text{~otherwise},\notag
\end{cases}
\end{align}
where $\lambda,\gamma,\delta$ are constants defined as in Equation \eqref{equ 1.1 generator}.

For any $t\geq 0$, we define $\widehat{C}_t=\{x:~\xi_t(x)=2\}$ and $\widehat{D}_t=\{x:~\xi_t(x)=1\}$. We write $\xi_t, \widehat{C}_t, \widehat{D}_t$ as $\xi_t^{(C,D)}, \widehat{C}_t^{(C,D)}, \widehat{D}_t^{(C,D)}$ when
$\widehat{C}_0=C$ and $\widehat{D}_0=D$.

The following proposition gives a duality relationship between the two-stage contact process $\{\eta_t\}_{t\geq 0}$ and the `on-off' process $\{\xi_t\}_{t\geq 0}$ with identical parameters $\lambda,\delta, \gamma$, which was proved by Krone in \cite{Krone1999}.
\begin{proposition}\label{proposition 5.2 1999}
(Krone, 1999) For any $A,B,C,D\subseteq \mathbb{Z}^d$ that $A\bigcap B=\emptyset$ and $C\bigcap D=\emptyset$,
\begin{align*}
&P_\lambda\Big(\eta_t^{(C,D)}(x)=2 \text{~for some~} x\in A\text{~or~}\eta_t^{(C,D)}(y)\neq 0\text{~for some~} y\in B\Big)\\
&=P_\lambda\Big(\xi_t^{(B,A)}(x)=2 \text{~for some~} x\in D\text{~or~}\xi_t^{(B,A)}(y)\neq 0\text{~for some~} y\in C\Big).
\end{align*}
\end{proposition}
If we let $C=\mathbb{Z}^d$ and $D=\emptyset$ while let $t$ grow to infinity, then we have the following direct corollary.
\begin{corollary}\label{corollary 5.3 1999}
(Krone, 1999) For any $A,B\subseteq \mathbb{Z}^d$ that $A\bigcap B=\emptyset$,
\begin{align*}
1-\pi(A,B,\lambda,d)=P_\lambda\Big(\widehat{C}_t^{(B,A)}\bigcup \widehat{D}_t^{(B,A)}\neq \emptyset\text{~for all~}t\geq 0\Big).
\end{align*}
\end{corollary}
To bound $P_\lambda\Big(\widehat{C}_t^{(B,A)}\bigcup \widehat{D}_t^{(B,A)}\neq \emptyset\text{~for all~}t\geq 0\Big)$ from above, we introduce a two-type branching process where there are some type $1$ individuals and some type $2$ individuals at $t=0$. Each individual is independently removed from the system at rate $1$. Each type $1$ individual independently becomes a type $2$ individual at rate $\gamma$. Each type $2$ individual independently becomes a type $1$ individual at rate $\delta$ while gives birth to a type $1$ individual at rate $\lambda$.

That is to say, if we use $\widehat{\zeta}_t$ to denote the number of type $2$ individuals at $t$ while use $\widehat{g}_t$ to denote the number of type $1$ individuals at $t$, then $\{\big(\widehat{\zeta}_t,\widehat{g}_t\big)\}_{t\geq 0}$ evolves as follows.
\begin{equation}\label{equ 5.3}
\big(\widehat{\zeta}_t,\widehat{g}_t\big)\text{~flips to~}
\begin{cases}
\big(\widehat{\zeta}_t-1,\widehat{g}_t\big) &\text{~at rate~}\widehat{\zeta}_t,\\
\big(\widehat{\zeta}_t,\widehat{g}_t-1\big) &\text{~at rate~}\widehat{g}_t,\\
\big(\widehat{\zeta}_t+1,\widehat{g}_t-1\big) &\text{~at rate~}\gamma\widehat{g}_t,\\
\big(\widehat{\zeta}_t-1,\widehat{g}_t+1\big) &\text{~at rate~}\delta\widehat{\zeta}_t,\\
\big(\widehat{\zeta}_t,\widehat{g}_t+1\big) &\text{~at rate~}\lambda\widehat{\zeta}_t,\\
0 &\text{~otherwise}.
\end{cases}
\end{equation}

For $m, n\geq 0$, we use $\widehat{\pi}(n,m)$ to denote the probability of the event that $\widehat{\zeta}_t+\widehat{g}_t>0$ for all $t\geq 0$ conditioned on there being $n$ type $2$ individuals and $m$ type $1$ individuals at $t=0$. We write $\widehat{\pi}(n,m)$ as $\widehat{\pi}(n,m,\lambda)$ when we need to point out the rate $\lambda$ at which a type $2$ individual gives birth to a type $1$ individual. Then, we have the following lemma.

\begin{lemma}\label{lemma 5.2}
For any $m,n\geq 0$ and $\lambda>\frac{1+\delta+\gamma}{\gamma}$,
\[
\widehat{\pi}(n,m,\lambda)=1-
\big(\frac{1+\delta+\gamma}{\lambda\gamma}\big)^n\big(1-\frac{\lambda\gamma-(1+\delta+\gamma)}{\lambda(\gamma+1)}\big)^m.
\]
\end{lemma}

\proof

According to the property of independent exponential times and the strong Markov property,
\begin{equation}\label{equ 5.4}
\begin{cases}
\widehat{\pi}(1,0)&=\frac{\lambda}{1+\delta+\lambda}\widehat{\pi}(1,1)+\frac{\delta}{1+\delta+\lambda}\widehat{\pi}(0,1),\\
\widehat{\pi}(0,1)&=\frac{\gamma}{1+\gamma}\widehat{\pi}(1,0).
\end{cases}
\end{equation}
Since the activities of different individuals are independent, for any $m,n\geq 0$,
\begin{equation}\label{equ 5.5}
\widehat{\pi}(n,m)=1-\big(1-\widehat{\pi}(1,0)\big)^n\big(1-\widehat{\pi}(0,1)\big)^m.
\end{equation}
Applying Equations \eqref{equ 5.4} and \eqref{equ 5.5} with $m=n=1$, we have
\begin{equation}\label{equ 5.6}
\widehat{\pi}(1,0)\big[\lambda\gamma\big(1-\widehat{\pi}(1,0)\big)-(1+\delta+\gamma)\big]=0.
\end{equation}

By direct calculation, when $\lambda>\frac{1+\delta+\gamma}{\gamma}$, the mean of the number of type $2$ children of a type $2$ father
is
\[
\frac{\lambda\gamma}{1+\delta+\gamma}>1.
\]
Therefore, $\widehat{\pi}(1,0,\lambda)>0$ when $\lambda>\frac{1+\delta+\gamma}{\gamma}$ according to the classic theory of branching processes. Then, by Equation \eqref{equ 5.6},
\[
\widehat{\pi}(1,0,\lambda)=1-\frac{1+\delta+\gamma}{\lambda\gamma}
\]
when $\lambda>\frac{1+\delta+\gamma}{\gamma}$. As a result, by Equations \eqref{equ 5.4} and \eqref{equ 5.5},
\[
\widehat{\pi}(0,1,\lambda)=\frac{\lambda\gamma-(1+\delta+\gamma)}{\lambda(\gamma+1)}
\]
and
\begin{equation}
\widehat{\pi}(n,m,\lambda)=1-\big(\frac{1+\delta+\gamma}{\lambda\gamma}\big)^n\big(1-\frac{\lambda\gamma-(1+\delta+\gamma)}{\lambda(\gamma+1)}\big)^m
\end{equation}
for any $m,n\geq 0$ and $\lambda>\frac{1+\delta+\gamma}{\gamma}$.

\qed

Now we can give the proof of Lemma \ref{lemma 5.1}.

\proof[Proof of Lemma \ref{lemma 5.1}]

For the `on-off' process $\{\xi_t\}_{t\geq 0}$ on $\mathbb{Z}^d$ with parameter $\frac{\lambda}{2d},\delta,\gamma$, a type $2$ vertex gives birth to a type $1$ vertex at rate
\[
\frac{\lambda}{2d}\sum_{y\sim x}1_{\{\xi_t(y)=0\}}\leq \frac{\lambda}{2d}\times 2d=\lambda.
\]
As a result, for $A, B\subseteq \mathbb{Z}^d$ that $|A|=m$ and $|B|=n$ while $A\bigcap B=\emptyset$, $\widehat{C}_t^{(B,A)}$ and $\widehat{D}_t^{(B,A)}$ are stochastic denominated from above by the numbers of type $2$ individuals and type $1$ individuals at moment $t$ respectively of the two-type branching process with $n$ initial type $2$ individuals and $m$ initial type $1$ individuals. Therefore,
\begin{equation}\label{equ 5.7}
P_{\frac{\lambda}{2d},d}\Big(\widehat{C}_t^{(B,A)}\bigcup \widehat{D}_t^{(B,A)}\neq \emptyset\text{~for all~}t\geq 0\Big)\leq \widehat{\pi}(n,m,\lambda).
\end{equation}
Lemma \ref{lemma 5.1} follows from Corollary \ref{corollary 5.3 1999}, Lemma \ref{lemma 5.2} and Equation \eqref{equ 5.7} directly.

\qed

\section{Lower bounds of $\text{~}1-\pi(A,B)$}\label{section six}
In this section we will give lower bounds of $1-\pi(A,B)$ to accomplish the proof of Theorem \ref{theorem 2.1 new main}. First we introduce some notations and definitions for later use. Let $\widetilde{X}_1,\ldots,\widetilde{X}_n,\ldots$ be independent and identically distributed random variables that
\[
P(\widetilde{X}_1=1)=e^{-(1+\delta)}(1-e^{-\gamma})=1-P(\widetilde{X}_1=0),
\]
then, for each integer $M\geq 1$, we define
\[
\widetilde{\alpha}(M)=P\Big(\frac{\sum_{i=1}^M\widetilde{X}_i}{M}\geq \frac{e^{-(1+\delta)}(1-e^{-\gamma})}{2}\Big).
\]
For any $d\geq 1, \lambda>0$ and $n\geq 1$, we define
\begin{align*}
\widetilde{b}(d,n,\lambda)=\inf\Bigg\{&\nu_{\frac{\lambda}{2d},d}\Big(\eta(x)\neq 0\text{~for some~}x\in A\Big):\\
&A\subseteq \mathbb{Z}^d\text{~and~}|A|=n\Bigg\}.
\end{align*}
The aim of this section is to prove the following two lemmas.
\begin{lemma}\label{lemma 6.1}
For $\lambda>\frac{1+\gamma+\delta}{\gamma}$ and $n\geq 1$,
\[
\liminf_{d\rightarrow+\infty}\widetilde{b}(d,n,\lambda)\geq \frac{1}{\frac{1}{n}\frac{2(\gamma+1)}{\gamma-\frac{1+\delta+\gamma}{\lambda}}+\frac{n-1}{n}}.
\]
\end{lemma}

\begin{lemma}\label{lemma 6.2}
For $m,n\geq 0$, $\lambda>\frac{1+\delta+\gamma}{\gamma}$, $M>n+m$ and sufficiently large $d$,
\begin{align*}
1-\pi(A,B,\frac{\lambda}{2d},d)\geq & \Bigg\{1-\Bigg(\frac{1+\delta+\gamma}{\frac{(2d-M)\lambda}{2d}\gamma}\Bigg)^n\Bigg(1-\frac{\frac{(2d-M)\lambda}{2d}\gamma-(1+\delta+\gamma)}
{\frac{(2d-M)\lambda}{2d}(\gamma+1)}\Bigg)^m\Bigg\}\\
&\times \widetilde{\alpha}(M)\widetilde{b}\Big(d, \Big\lceil\frac{Me^{-(1+\delta)}(1-e^{-\gamma})}{2}\Big\rceil,\lambda\Big)
\end{align*}
for any $A, B\subseteq \mathbb{Z}^d$ that $|A|=m, |B|=n$ and $A\bigcap B=\emptyset$, where
\[
\lceil x\rceil=\inf\big\{m:~m\geq x\text{~and~}m\text{~is an integer}\big\}.
\]
\end{lemma}

Before proving Lemmas \ref{lemma 6.1} and \ref{lemma 6.2}, we first show how to utilize these two lemmas to prove Theorem \ref{theorem 2.1 new main}.

\proof[Proof of Theorem \ref{theorem 2.1 new main}]

For simplicity, we use $\widetilde{c}(M,d,\lambda)$ to denote
\begin{align*}
&\Bigg\{1-\Bigg(\frac{1+\delta+\gamma}{\frac{(2d-M)\lambda}{2d}\gamma}\Bigg)^n\Bigg(1-\frac{\frac{(2d-M)\lambda}{2d}\gamma-(1+\delta+\gamma)}
{\frac{(2d-M)\lambda}{2d}(\gamma+1)}\Bigg)^m\Bigg\}\\
&\times \widetilde{\alpha}(M)\widetilde{b}\Big(d, \Big\lceil\frac{Me^{-(1+\delta)}(1-e^{-\gamma})}{2}\Big\rceil,\lambda\Big)
\end{align*} while use $\mu(M)$ to denote $\Big\lceil\frac{Me^{-(1+\delta)}(1-e^{-\gamma})}{2}\Big\rceil$. Then, according to Lemmas \ref{lemma 5.1} and \ref{lemma 6.2},
\begin{equation}\label{equ 6.1}
\Pi(m,n,\lambda,d)\leq 1-\bigg(1-\frac{\lambda\gamma-(1+\delta+\gamma)}{\lambda(\gamma+1)}\bigg)^m\bigg(\frac{1+\delta+\gamma}{\lambda\gamma}\bigg)^n
-\widetilde{c}(M,d,\lambda)
\end{equation}
for $m,n\geq 0$ and $\lambda>\frac{1+\delta+\gamma}{\gamma}$. By Lemma \ref{lemma 6.1},
\begin{align}\label{equ 6.2}
\liminf_{d\rightarrow+\infty}\widetilde{c}(M,d,\lambda)\geq & \Bigg\{1-\bigg(1-\frac{\lambda\gamma-(1+\delta+\gamma)}{\lambda(\gamma+1)}\bigg)^m\bigg(\frac{1+\delta+\gamma}{\lambda\gamma}\bigg)^n\Bigg\} \\
&\times \widetilde{\alpha}(M)\frac{1}{\frac{2(\gamma+1)}{\mu(M)(\gamma-\frac{1+\delta+\gamma}{\lambda})}+\frac{\mu(M)-1}{\mu(M)}} \notag
\end{align}
for sufficiently large $M$ that $\mu(M)>1$ and $\lambda>\frac{1+\delta+\gamma}{\gamma}$. By Equations \eqref{equ 6.1} and \eqref{equ 6.2},
\begin{align}\label{equ 6.3}
\limsup_{d\rightarrow+\infty}\Pi(m,n,\lambda,d)\leq &
\Bigg\{1-\bigg(1-\frac{\lambda\gamma-(1+\delta+\gamma)}{\lambda(\gamma+1)}\bigg)^m\bigg(\frac{1+\delta+\gamma}{\lambda\gamma}\bigg)^n\Bigg\} \\
&\times \Bigg(1-\widetilde{\alpha}(M)\frac{1}{\frac{2(\gamma+1)}{\mu(M)(\gamma-\frac{1+\delta+\gamma}{\lambda})}
+\frac{\mu(M)-1}{\mu(M)}}\Bigg) \notag
\end{align}
for any sufficiently large $M$ and $\lambda>\frac{1+\delta+\gamma}{\gamma}$. According to the law of large numbers,
\[
\lim_{M\rightarrow+\infty}\widetilde{\alpha}(M)=1.
\]
As a result, let $M\rightarrow+\infty$,
\begin{equation}\label{equ 6.4}
\limsup_{d\rightarrow+\infty}\Pi(m,n,\lambda,d)\leq 0
\end{equation}
for any $\lambda>\frac{1+\delta+\gamma}{\gamma}$ according to Equation \eqref{equ 6.3} and the fact that
\[
\lim_{M\rightarrow+\infty}\mu(M)=+\infty.
\]
Since $\Pi(m,n,\lambda,d)$ is nonnegative, Theorem \ref{theorem 2.1 new main} follows from Equation \eqref{equ 6.4} directly.

\qed

Now we give the proof of Lemma \ref{lemma 6.1}.

\proof[Proof of Lemma \ref{lemma 6.1}]

Let $A\subseteq \mathbb{Z}^d$ that $|A|=n$, then according to the definition of $\nu$ and Cauchy-Schwartz's inequality,
\begin{align}\label{equ 6.5}
&\nu_{\frac{\lambda}{2d},d}\Big(\eta(x)\neq 0\text{~for some~}x\in A\Big) \\
&=\lim_{t\rightarrow+\infty}P_{\frac{\lambda}{2d},d}\Big(\eta_t^{(\mathbb{Z}^d,\emptyset)}(x)\neq 0\text{~for some~}x\in A\Big) \notag\\
&=\lim_{t\rightarrow+\infty}P_{\frac{\lambda}{2d},d}\Big(\zeta_t(x)+g_t(x)>0\text{~for some~}x\in A\Big) \notag\\
&\geq \limsup_{t\rightarrow+\infty}P_{\frac{\lambda}{2d},d}\Big(\zeta_t(x)>0\text{~for some~}x\in A\Big) \notag\\
&=\limsup_{t\rightarrow+\infty}P_{\frac{\lambda}{2d},d}\Big(\sum_{x\in A}\zeta_t(x)>0\Big) \geq \limsup_{t\rightarrow+\infty}\frac{\Big(E_{\frac{\lambda}{2d},d}\sum_{x\in A}\zeta_t(x)\Big)^2}{E_{\frac{\lambda}{2d},d}\Bigg(\Big(\sum_{x\in A}\zeta_t(x)\Big)^2\Bigg)}, \notag
\end{align}
where $\{\big(\zeta_t(x),g_t(x)\big):~t\geq 0, x\in \mathbb{Z}^d\}$ is our auxiliary model defined as in Section \ref{section four}. We have shown in Section \ref{section four} that $E\zeta_t(x)=E\zeta_t(O)\equiv1$, then by Equation \eqref{equ 6.5},
\begin{align}\label{equ 6.6}
&\nu_{\frac{\lambda}{2d},d}\Big(\eta(x)\neq 0\text{~for some~}x\in A\Big)\geq \limsup_{t\rightarrow+\infty}\frac{1}{\frac{1}{n^2}\sum_{x,y\in A}F_t(y-x,1)}\\
&=\limsup_{t\rightarrow+\infty}\frac{1}{\frac{1}{n}F_t(O,1)+\frac{1}{n^2}\sum_{x,y\in A, x\neq y}F_t(y-x,1)}, \notag
\end{align}
where $F_t(x,1)=E_{\frac{\lambda}{2d},d}\Big(\zeta_t(O)\zeta_t(x)\Big)=E_{\frac{\lambda}{2d},d}\Big(\zeta_t(y)\zeta_t(x+y)\Big)$ defined as in Section \ref{section four}. For $\lambda>\frac{1+\delta+\gamma}{\gamma}$ and sufficiently large $d$, let $K_{\frac{\lambda}{2d}}$ be the function on $X_4$ defined as before Equation \eqref{equ 4.19}, then $K_{\frac{\lambda}{2d}}$ is the eigenvector with respect to the eigenvalue $0$ of the $X_4\times X_4$ matrix $G_\frac{\lambda}{2d}$ and
\begin{equation}\label{equ 6.6 two}
F_t(x,1)\leq \frac{K_{\frac{\lambda}{2d}}(x,1)}{\inf_{(x,i)\in X_4} K_{\frac{\lambda}{2d}}(x,i)}
\end{equation}
as we have shown in Equation \eqref{equ 4.13} and the proof of Equation \eqref{equ 2.5}. Note that $\inf_{(x,i)\in X_4} K_{\frac{\lambda}{2d}}(x,i)>0$ when $\lambda>\frac{1+\gamma+\delta}{\gamma}$ and
$d$ is sufficiently large according to the definition of $K_\frac{\lambda}{2d}$ and the fact that
\[
\Gamma(e_1,1)\leq \widetilde{\Gamma}(e_1)=\frac{1}{2d}+O(\frac{1}{d^2}).
\]
As we have defined in Section \ref{section four},
\[
K_{\frac{\lambda}{2d}}(x,1)=\Gamma(x,1)+h_{\frac{\lambda}{2d}},
\]
where
\[
h_{\frac{\lambda}{2d}}=\frac{\gamma[1-2\Gamma(O,2)]-2\Gamma(e_1,2)-\frac{1+\delta+\gamma}{\lambda}}{\gamma+2+\frac{1+\delta+\gamma}{\lambda}}.
\]
By Equation \eqref{equ 4.19},
\[
\inf_{(x,i)\in X_4}K_{\frac{\lambda}{2d}}(x,i)\geq \inf\big\{h_{\frac{\lambda}{2d}}, \gamma[1-2\Gamma(e_1,1)-h_{\frac{\lambda}{2d}}]\big\}.
\]
Then, according to the definition of $h_\lambda$ and the fact that $\Gamma(e_1,1)\leq \widetilde{\Gamma}(e_1)=\frac{1}{2d}+O(\frac{1}{d^2})$,
\[
\gamma[1-2\Gamma(e_1,1)-h_{\frac{\lambda}{2d}}]>h_{\frac{\lambda}{2d}}
\]
and hence
\begin{equation}\label{equ 6.7}
\inf_{(x,i)\in X_4}K_{\frac{\lambda}{2d}}(x,i)\geq h_{\frac{\lambda}{2d}}
\end{equation}
for sufficiently large $d$. By Equations \eqref{equ 6.6 two} and \eqref{equ 6.7}, for sufficiently large $d$,
\begin{align}\label{equ 6.8}
F_t(x,1)&\leq \frac{\Gamma(x,1)+h_{\frac{\lambda}{2d}}}{h_{\frac{\lambda}{2d}}} \leq
\frac{\widetilde{\Gamma}(x)+h_{\frac{\lambda}{2d}}}{h_{\frac{\lambda}{2d}}} \leq \frac{\widetilde{\Gamma}(e_1)+h_{\frac{\lambda}{2d}}}{h_{\frac{\lambda}{2d}}}
\end{align}
for any $x\neq O$ while
\begin{equation}\label{equ 6.10}
F_t(O,1)\leq \frac{1+h_{\frac{\lambda}{2d}}}{h_{\frac{\lambda}{2d}}}.
\end{equation}
By Equations \eqref{equ 6.6}, \eqref{equ 6.8} and \eqref{equ 6.10},
\begin{equation}\label{equ 6.11}
\nu_{\frac{\lambda}{2d},d}\Big(\eta(x)\neq 0\text{~for some~}x\in A\Big)\geq \frac{1}{\frac{1}{n}\frac{1+h_{\frac{\lambda}{2d}}}{h_{\frac{\lambda}{2d}}}
+\frac{n(n-1)}{n^2}\frac{\widetilde{\Gamma}(e_1)+h_{\frac{\lambda}{2d}}}{h_{\frac{\lambda}{2d}}}}.
\end{equation}
Since $\widetilde{\Gamma}(e_1)=\frac{1}{2d}+O(\frac{1}{d^2})$ and
\[
\lim_{d\rightarrow+\infty}h_{\frac{\lambda}{2d}}=\frac{\gamma-\frac{1+\delta+\gamma}{\lambda}}{\gamma+2+\frac{1+\delta+\gamma}{\lambda}},
\]
Lemma \ref{lemma 6.1} follows directly from Equation \eqref{equ 6.11}.

\qed

At last, we give the proof of Lemma \ref{lemma 6.2}.

\proof[Proof of Lemma \ref{lemma 6.2}]
Throughout this proof we assume that $n,m,M,d$ are fixed such that $M>n+m$ and $2d>M$.
For $A,B\subseteq \mathbb{Z}^d$ such that $|A|=m, |B|=n$ and $A\bigcap B=\emptyset$, let
\[
\tau_M(A,B)=\inf\{t\geq 0:~|\widehat{C}_t^{(B,A)}|+|\widehat{D}_t^{(B,A)}|=M\},
\]
where $\{\xi_t\}_{t\geq 0}$ is the `on-off' process introduced in Section \ref{section five} and
\[
\widehat{C}_t=\{x:~\xi_t(x)=2\}\text{~while~} \widehat{D}_t=\{x:~\xi_t(x)=1\}
\]
defined as in Section \ref{section five}. Let $\{\big(\widehat{\zeta}_t^M, \widehat{g}_t^M\big)\}_{t\geq 0}$ be the two-type branching process defined as in Section \ref{section five} with parameter $\frac{2d-M}{2d}\lambda,\gamma,\delta$, then we define
\[
\widehat{\tau}_M=\inf\{t\geq 0:~\widehat{\zeta}_t^M+\widehat{g}_t^M=M\}.
\]
For the `on-off' process $\{\xi_t^{(B,A)}\}_{t\geq 0}$ on $\mathbb{Z}^d$ with parameter $\frac{\lambda}{2d}$, a type $2$ vertex gives birth to a type $1$ vertex at rate at least
\[
\frac{\lambda}{2d}\times(2d-M)
\]
before the moment $\tau_M(A,B)$, since there are at lest $(2d-M)$ neighbors in state $0$ before the moment $\tau_M(A,B)$. As a result, $|\widehat{C}_t^{(B,A)}|+|\widehat{D}_t^{(B,A)}|$ is stochastic dominated from below by
$\widehat{\zeta}_t^M+\widehat{g}_t^M$ with $\widehat{\zeta}_0^M=n$ and $\widehat{g}_0^M=m$ for $t\in [0,\tau_M(A,B))$. Therefore,
\begin{equation}\label{equ 6.12}
P_{\frac{\lambda}{2d},d}\Big(\tau_M(A,B)<+\infty\Big)\geq P\Big(\widehat{\tau}_M<+\infty\Big|\widehat{\zeta}_0^M=n, \widehat{g}_0^M=m\Big).
\end{equation}
We claim that
\begin{align}\label{equ 6.13}
&P\Big(\widehat{\tau}_M<+\infty\Big|\widehat{\zeta}_0^M=n, \widehat{g}_0^M=m\Big)\geq \\
&P\Big(\widehat{\zeta}_t^M+\widehat{g}_t^M>0\text{~for all~}t\geq 0\Big|\widehat{\zeta}_0^M=n, \widehat{g}_0^M=m\Big). \notag
\end{align}
Equation \eqref{equ 6.13} holds according to the following analysis. Let
\[
\widetilde{\Xi}=\inf\Bigg\{P\Big(\widehat{\zeta}_{\frac{1}{2}}^M+\widehat{g}_{\frac{1}{2}}^M=0\Big|\widehat{\zeta}_0^M=k, \widehat{g}_0^M=l\Big):~l+k\leq M\Bigg\},
\]
then $\widetilde{\Xi}>0$ since $P\Big(\widehat{\zeta}_{\frac{1}{2}}^M+\widehat{g}_{\frac{1}{2}}^M=0\Big|\widehat{\zeta}_0^M=k, \widehat{g}_0^M=l\Big)>0$ for each pair of $(l,k)$ and there are finite many pairs of $(l,k)$s satisfying $l+k\leq M$. Let
\[
\widetilde{\tau}=\inf\big\{n:~n\text{~is a nonnegative integer and~}\widehat{\zeta}_{n+\frac{1}{2}}^M+\widehat{g}_{n+\frac{1}{2}}^M=0\big\},
\]
then $\widehat{\zeta}_{n}^M+\widehat{g}_{n}^M\leq M$ for each integer $n\geq 0$ on the event $\{\widehat{\tau}_M=+\infty\}$ and hence $\widetilde{\tau}$ is stochastic dominated from above by the random variable $\widetilde{Y}$ satisfying
\[
P(\widetilde{Y}=n)=\widetilde{\Xi}(1-\widetilde{\Xi})^n
\]
for $n=0,1,2,\ldots$ on the event $\{\widehat{\tau}_M=+\infty\}$. As a result,
\[
P\Big(\widetilde{\tau}<+\infty\Big|\widehat{\tau}_M=+\infty\Big)\geq P(\widetilde{Y}<+\infty)=1.
\]
That is to say
\[
\Big\{\widehat{\zeta}_t^M+\widehat{g}_t^M=0\text{~for some~}t\geq 0\Big\}\supseteq \Big\{\widehat{\tau}_M=+\infty\Big\}
\]
in the sense of ignoring a set with probability zero, Equation \eqref{equ 6.13} follows from which directly.

By Lemma \ref{lemma 5.2} with $\lambda$ replaced by $\frac{2d-M}{2d}\lambda$,
\begin{align*}
&P\Big(\widehat{\zeta}_t^M+\widehat{g}_t^M>0\text{~for all~}t\geq 0\Big|\widehat{\zeta}_0^M=n, \widehat{g}_0^M=m\Big)\\
&=1-\Bigg(\frac{1+\delta+\gamma}{\frac{(2d-M)\lambda}{2d}\gamma}\Bigg)^n\Bigg(1-\frac{\frac{(2d-M)\lambda}{2d}\gamma-(1+\delta+\gamma)}
{\frac{(2d-M)\lambda}{2d}(\gamma+1)}\Bigg)^m.
\end{align*}
Then, by Equations \eqref{equ 6.12} and \eqref{equ 6.13},
\begin{align}\label{equ 6.14}
&P_{\frac{\lambda}{2d},d}\Big(\tau_M(A,B)<+\infty\Big) \\
&\geq 1-\Bigg(\frac{1+\delta+\gamma}{\frac{(2d-M)\lambda}{2d}\gamma}\Bigg)^n\Bigg(1-\frac{\frac{(2d-M)\lambda}{2d}\gamma-(1+\delta+\gamma)}
{\frac{(2d-M)\lambda}{2d}(\gamma+1)}\Bigg)^m.\notag
\end{align}
According to a similar analysis with that of the two-stage contact process given in \cite{Krone1999}, the `on-off' process is also monotonic with respect to the partial order `$\preceq$' defined in Section \ref{section two}. As a result, for $A_1, B_1\subseteq \mathbb{Z}^d$ that $A_1\bigcap B_1=\emptyset$ while $|A_1|+|B_1|=M$,
\begin{align}\label{equ 6.15}
&P_{\frac{\lambda}{2d},d}\Big(\widehat{C}_t^{(B_1,A_1)}\bigcup \widehat{D}_t^{(B_1,A_1)}\neq \emptyset\text{~for all~}t\geq 0\Big)\geq \\
&P_{\frac{\lambda}{2d},d}\Big(\widehat{C}_t^{(\emptyset,A_1\bigcup B_1)}\bigcup \widehat{D}_t^{(\emptyset,A_1\bigcup B_1)}\neq \emptyset\text{~for all~}t\geq 0\Big). \notag
\end{align}
By direct calculation, an initial type $1$ vertex becomes a type $2$ vertex at some moment $s<1$ and then stays in state $2$ till moment $t=1$ with probability at least
\[
P(\widetilde{Y}_2>1,\widetilde{Y}_3<1,\widetilde{Y}_4>1)=e^{-(1+\delta)}(1-e^{-\gamma}),
\]
where $\widetilde{Y}_2, \widetilde{Y}_3, \widetilde{Y}_4$ are independent exponential times with rates $1, \gamma, \delta$ respectively. Therefore, for $A_1, B_1\subseteq \mathbb{Z}^d$ that $A_1\bigcap B_1=\emptyset$ while $|A_1|+|B_1|=M$, $|\widehat{C}_1^{(\emptyset,A_1\bigcup B_1)}|$ is stochastic dominated from below by a random variable following the binomial distribution $B(M,e^{-(1+\delta)}(1-e^{-\gamma}))$ and
\begin{equation}\label{equ 6.16}
P_{\frac{\lambda}{2d},d}\Big(|\widehat{C}_1^{(\emptyset,A_1\bigcup B_1)}|\geq \frac{Me^{-(1+\delta)}(1-e^{-\gamma})}{2} \Big)\geq \widetilde{\alpha}(M).
\end{equation}
For $A_2, B_2\subseteq \mathbb{Z}^d$ that $|B_2|=\big\lceil\frac{Me^{-(1+\delta)}(1-e^{-\gamma})}{2}\big\rceil$ and $A_2\bigcap B_2=\emptyset$, by Corollary \ref{corollary 5.3 1999} and the definition of $\widetilde{b}(d, \big\lceil\frac{Me^{-(1+\delta)}(1-e^{-\gamma})}{2}\big\rceil, \lambda)$,
\begin{align}\label{equ 6.17}
&P_{\frac{\lambda}{2d},d}\Big(\widehat{C}_t^{(B_2,A_2)}\bigcup \widehat{D}_t^{(B_2,A_2)}\neq \emptyset\text{~for all~}t\geq 0\Big)   \\
&\geq P_{\frac{\lambda}{2d},d}\Big(\widehat{C}_t^{(B_2,\emptyset)}\bigcup \widehat{D}_t^{(B_2,\emptyset)}\neq \emptyset\text{~for all~}t\geq 0\Big) \notag\\
&=1-\pi(\emptyset, B_2,\frac{\lambda}{2d},d)=\nu_{\frac{\lambda}{2d},d}\Big(\eta(x)\neq 0\text{~for some~}x\in B_2\Big) \notag\\
&\geq \widetilde{b}(d, \big\lceil\frac{Me^{-(1+\delta)}(1-e^{-\gamma})}{2}\big\rceil,\lambda). \notag
\end{align}
By Equations \eqref{equ 6.16}, \eqref{equ 6.17} and the Markov property, for $A_1, B_1\subseteq \mathbb{Z}^d$ that $A_1\bigcap B_1=\emptyset$ while $|A_1|+|B_1|=M$,
\begin{align}\label{equ 6.18}
&P_{\frac{\lambda}{2d},d}\Big(\widehat{C}_t^{(B_1,A_1)}\bigcup \widehat{D}_t^{(B_1,A_1)}\neq \emptyset\text{~for all~}t\geq 0\Big) \\
&\geq \widetilde{\alpha}(M)\widetilde{b}\Big(d, \Big\lceil\frac{Me^{-(1+\delta)}(1-e^{-\gamma})}{2}\Big\rceil,\lambda\Big). \notag
\end{align}
On the event $\tau_M(A,B)<+\infty$, $\widetilde{B}:=\widehat{C}_{\tau_M(A,B)}^{(B,A)}$ and $\widetilde{A}:=\widehat{D}_{\tau_M(A,B)}^{(B,A)}$ satisfy $\widetilde{A}\bigcap \widetilde{B}=\emptyset$ while $|\widetilde{A}|+|\widetilde{B}|=M$. Therefore, by Equations \eqref{equ 6.14}, \eqref{equ 6.18} and the strong Markov property,
\begin{align}\label{equ 6.19}
&P_{\frac{\lambda}{2d},d}\Big(\widehat{C}_t^{(B,A)}\bigcup \widehat{D}_t^{(B,A)}\neq \emptyset\text{~for all~}t\geq 0\Big) \\
&\geq \Bigg\{1-\Bigg(\frac{1+\delta+\gamma}{\frac{(2d-M)\lambda}{2d}\gamma}\Bigg)^n\Bigg(1-\frac{\frac{(2d-M)\lambda}{2d}\gamma-(1+\delta+\gamma)}
{\frac{(2d-M)\lambda}{2d}(\gamma+1)}\Bigg)^m\Bigg\} \notag\\
&\times \widetilde{\alpha}(M)\widetilde{b}\Big(d, \Big\lceil\frac{Me^{-(1+\delta)}(1-e^{-\gamma})}{2}\Big\rceil,\lambda\Big) \notag
\end{align}
for $A,B\subseteq \mathbb{Z}^d$ that $|A|=m, |B|=n$ and $A\bigcap B=\emptyset$. Lemma \ref{lemma 6.2} follows from Corollary \ref{corollary 5.3 1999} and Equation \eqref{equ 6.19} directly.

\qed

\quad

\textbf{Acknowledgments.}
The author is grateful to the financial
support from the National Natural Science Foundation of China with
grant number 11501542 and the financial support from Beijing
Jiaotong University with grant number KSRC16006536.

{}
\end{document}